%% file: stochasticcoildesign.tex
\documentclass[10pt]{iopart}
\usepackage[citestyle=alphabetic,bibstyle=authortitle]{biblatex}
\expandafter\let\csname equation*\endcsname\relax
\expandafter\let\csname endequation*\endcsname\relax
\usepackage{mathtools}
\usepackage{bbm}
\usepackage{iopams}  
\include{latex-includes/makros}
\usepackage{xcolor}
\definecolor{blue}{rgb}{0.2980392156862745, 0.4470588235294118, 0.6901960784313725}
\definecolor{green}{rgb}{0.3333333333333333, 0.6588235294117647, 0.40784313725490196}
\definecolor{red}{rgb}{0.7686274509803922, 0.3058823529411765, 0.3215686274509804}
\definecolor{purple}{rgb}{0.5058823529411764, 0.4470588235294118, 0.6980392156862745}
\definecolor{yellow}{rgb}{0.8, 0.7254901960784313, 0.4549019607843137}
\definecolor{lightblue}{rgb}{0.39215686274509803, 0.7098039215686275, 0.803921568627451}
\definecolor{orange}{rgb}{0.8666666666666667, 0.5176470588235295, 0.3215686274509804}
\definecolor{brown}{rgb}{0.5764705882352941, 0.47058823529411764, 0.3764705882352941}
\definecolor{pink}{rgb}{0.8549019607843137, 0.5450980392156862, 0.7647058823529411}
\definecolor{gray}{rgb}{0.5490196078431373, 0.5490196078431373, 0.5490196078431373}

\usepackage{float}
\usepackage{standalone}
\usepackage{pgfplots}
\usepackage{pgfplotstable}
\usepgfplotslibrary{groupplots}
\usetikzlibrary{pgfplots.groupplots}
\usetikzlibrary{calc,arrows,spy,intersections}
\pgfplotsset{compat=1.13}
\begin{document}
\providecommand{\datadir}{plots/}%

\title[Stochastic Coil Design]{Single-stage gradient-based stellarator coil design: stochastic optimization}

\author{Florian Wechsung$^1$, Andrew Giuliani$^1$, Matt Landreman$^2$, Antoine Cerfon$^1$, Georg Stadler$^1$}

\address{$^1$ Courant Institute of Mathematical Sciences, New York University, New York, USA}
\address{$^2$ University of Maryland-College Park, Maryland, USA}
\ead{wechsung@nyu.edu, giuliani@cims.nyu.edu, mattland@umd.edu, cerfon@cims.nyu.edu, stadler@cims.nyu.edu}
\vspace{10pt}
\begin{indented}
\item[]\today
\end{indented}

\begin{abstract}
  We extend the single-stage stellarator coil design approach for quasi-symmetry on axis from [Giuliani et al, 2020] to additionally take into account coil manufacturing errors. By modeling coil errors independently from the coil discretization,
  we have the flexibility to consider realistic forms of coil errors. The corresponding stochastic optimization problems are formulated using a risk-neutral approach
  and risk-averse approaches. We present an efficient, gradient-based descent algorithm which relies on analytical derivatives to solve these problems. In a
  comprehensive numerical study, we compare the coil designs resulting from deterministic and risk-neutral stochastic optimization and find that the risk-neutral formulation results in more robust configurations and reduces the number of local minima of the optimization problem. We also compare deterministic and risk-neutral approaches in terms of quasi-symmetry on and away from the magnetic axis, and in terms of the confinement of particles released close to the axis. Finally, we show that for the optimization problems we consider, a risk-averse objective using the Conditional Value-at-Risk leads to results which are similar to the risk-neutral objective.
\end{abstract}

%
%
%
%
%

\section{Introduction}

The design, manufacturing, and assembly of the primary coil system for stellarators are among the most technically challenging aspects of the construction of stellarators \cite{Strykowsky09, Neilson10,Klinger2013}, and represent a large fraction of the total construction cost. Some of these challenges are intrinsic to the nature of stellarators: non-axisymmetric coil systems are expected to be more complex than the axisymmetric coil systems of tokamaks \cite{Strykowsky09,gates2018stellarator}. However, some of these challenges are directly linked to the design optimization process, and the optima which have been selected for construction. Specifically, the emphasis has historically been given to optimization metrics corresponding to plasma performance, and less attention was given to the engineering requirements and constraints to achieve such performance. Stellarator designs considered optimal from a physics point of view could only be realized with complex and expensive coil configurations. The difficulties encountered in the construction of large scale stellarators \cite{Strykowsky09, Neilson10,Klinger2013} has recently triggered a renewed research effort toward the development of tools enabling the design of simpler and more efficient coil systems that still lead to strong plasma performance \cite{Landreman2016,Landreman_2017,Hudson2018,Paul_2018,Carlton2021}.

Beside the engineering complexity of a given design, the lack of
robustness of the physics performance in the presence of coil
manufacturing and assembly errors is also a strong driver for the
cost of a machine, since it requires tight tolerances at every step of
the manufacturing and assembly process. Both challenges are
related but not identical: a relatively simple coil system whose performance
degrades strongly with manufacturing and alignment errors may not
be as desirable as a more
complex coil system with more robust performance. Work has therefore also been lately devoted to the
design of efficient numerical methods for the evaluation of the
sensitivity of error fields to coil perturbations
\cite{Zhu_2018_error,Zhu_2019}, and the sensitivity of physical
quantities to error fields
\cite{Landreman_2018_sensitivity,geraldini_2021}. These methods can be
included in deterministic stellarator optimization codes, and serve to
narrow the search to configurations with lower sensitivity or with
sensitivity with respect to perturbations that are more easily
controlled. 

A complementary approach to deal with the challenge of
strong sensitivity and tight tolerances is to account for errors during the optimization process, via
stochastic optimization. These errors may either be engineering errors or possibly errors due to the limitations of
the physics models used. In stochastic optimization, the objective is a
function both of the controls and of the random model errors and hence is
itself a random variable.  This is
precisely the approach we adopt in this article, with the randomness
corresponding to perturbations of the location and shape of the coils. We observe that while stochastic optimization is not yet a widely applied method for the design of magnetic confinement devices, it was used in the design of the CNT stellarator ~\cite{Kremer_2003, Pedersen_2006}. CNT consists of four circular coils, namely two interlocking (IL) coils and two poloidal field (PF) coils, and the CNT optimization problem involves the angle between the IL coils, as well as the current ratio between the IL and the PF coils. Coil errors involving tilts and shifts of these coils were considered and the goal was to optimize the average volume with good flux surfaces.
More recently, for the problem of designing coils corresponding to a desired plasma boundary, Lobsien~\textit{et al.}~\cite{lobsien_physics_2020,lobsien_improved_2020} demonstrated that stochastic optimization leads to simpler and better performing coils.

In this manuscript we present a stochastic version of the single-stage coil design framework recently introduced in \cite{giuliani2020}. As the purpose of our work is mainly to introduce a new paradigm for stochastic coil design optimization, the physics basis for the design is simple: we consider a vacuum field, and  optimize directly for a target value of the rotational transform and for quasi-symmetry on the magnetic axis. The main contribution of this work is a formulation that enables us to consider more general and realistic coil perturbations than considered thus far, to implement gradient based optimisation algorithms relying on analytic derivatives, and to compare the performance of different forms of stochastic optimization.  

We find that both the stochastic and
deterministic formulations result in different
designs depending on the initialization of the optimization algorithm, which indicates the existence of multiple local minima. However, the variability of the designs
corresponding to the stochastic problem is substantially reduced as compared to the designs obtained from the deterministic optimization problem. The coil systems we obtain for different initial conditions are much more similar to one another when using stochastic optimization than with deterministic optimization. We compare the performance of the obtained configurations in the presence of coil errors by evaluating the level of quasi symmetry near and away from the axis, the rotational transform on axis as well as particle loss fraction. We observe that the configurations found by stochastic optimization outperform those obtained from the deterministic formulation. Furthermore, the different minima obtained from stochastic minimization all perform very similarly, which suggests that the optimization algorithm does not get trapped in poor local minima. For all these reasons, our work demonstrates the strong potential of stochastic optimization for stellarator design, and motivates its application to more detailed reactor design studies.



The structure of the article is as follows. In Section \ref{sec:model_perturb}, we present our mathematical description of the coils, and explain how we model random coil perturbations. In Section \ref{sec:stochastic}, we provide a brief summary of stochastic optimization, with a description of several variants which are relevant to stellarator design. We then review in Section \ref{sec:single_stage} the direct coil design paradigm that we modify for stochastic optimization, first introduced in~\cite{giuliani2020}. We present our main numerical results in Section~\ref{sec:num}, and summarize our work in Section~\ref{sec:conclusion}, where we also suggest directions for future work.

\section{Modeling coil perturbations}\label{sec:model_perturb}

\subsection{Physical representation of coil perturbations}
For this work, we make the common assumption in coil design that we can represent the coils as current-carrying filaments, i.e., we neglect the non-zero thickness of the coils, and simply model coils as curves in space.
A coil is then described by a periodic function $\Gammab : [0, 2\pi)\to \re^3$.
A standard approach to discretizing such coils is given by a truncated Fourier expansion,
that is the $j$-th coordinate of the $i$-th coil $\Gammab^{(i)}$ is given by
\begin{equation}
    \Gamma_j^{(i)}(\theta) = c_{j, 0}^{(i)} + \sum_{l=1}^{n_p^\text{coil}} s_{j, l}^{(i)} \sin(l\theta) + \sum_{l=1}^{n_p^\text{coil}} c_{j, l}^{(i)} \cos(l \theta).
\end{equation}
We collect the degrees of freedom for coil $i$ in the vector $\cb^{(i)}\in \re^{3 (2n_p^{\text{coils}} + 1)}$.
This approach is also used in the \pkg{FOCUS} code~\cite{zhu2017}.
We note that this formulation allows the coils to move freely in space, as compared to coil optimization codes that restrict the coils to lie on the so-called \emph{winding surface}.
This latter approach is employed in the \pkg{ONSET}~\cite{drevlak_optimization_1999}, \pkg{COILOPT}~\cite{strickler_designing_2002}, \pkg{COILOPT++}~\cite{gates_recent_2017} codes.

A straightforward way of modelling errors is to perturb the vector containing the degrees of freedom, $\tilde\cb^{(i)} = \cb^{(i)} + \epsb$ where $\epsb$ is a vector of independent random variables.
In the context of coil optimization, this approach has been used in~\cite{lobsien_stellarator_2018, lobsien_physics_2020, lobsien_improved_2020} which builds on \pkg{ONSET} and uses splines to represent coils.
There the spline anchor points were perturbed by independent, centered Gaussian random variables with small variance.
This approach does not require any modification of the objective
function implementation and hence a deterministic code can be extended towards
stochastic optimization with little effort.

However, coil errors originate from the manufacturing process and are thus independent of the coil description used in the design process.
Thus, it is unlikely that manufacturing errors satisfy stellarator symmetry~\cite{Dewar1998}, a property that most optimization studies and coil design codes assume for simplicity~\cite{Drevlak_2013,zhu2017,bader_2019,Henneberg_2019}.
Hence, using the same parametrization
for coil errors as used for the coils lacks generality, and may be unphysical. For instance, when Fourier modes are used for coils \emph{and} coil errors, these
errors affect coils globally and repeat themselves along the coils,
with more repetitions with increasing Fourier mode numbers.  Using
splines to describe the coil geometry as well as manufacturing errors
has the advantage of allowing the description of local manufacturing
errors. However, if the number of spline anchor points changes,
the characteristics of the manufacturing error changes as well.
%
Thus, when
the description of the coils and the errors is the same, changing the
coil discretization also means changing the type of errors
considered and convergence to a limit when refining the
coil description is unclear.

For this reason, we separate  the discretization of the coils and the
modelling of perturbations by considering additive perturbations
modelled by Gaussian processes \cite{RasmussenWilliams06}. To the best
of our knowledge, this is a novel approach for stellarator coil
optimization applications. However, similar approaches have been used
in other areas. For instance, such an approach has been used in the context of airfoil
optimization in two dimensions
\cite{chen_new_2011,wang_conditional_2011,liu_quantification_2017}.

\subsection{Randomly perturbed coils via stochastic processes}
We model the perturbations of the coil by random, periodic functions $\Xib:[0,2\pi) \to \re^3$, and denote the perturbed coil by $\tilde\Gammab(\theta) = \Gammab(\theta) + \Xib(\theta)$.
We choose to model the components $(\Xi_1, \Xi_2, \Xi_3)$ of $\Xib$ as centered Gaussian
processes.

We briefly recall the definition and some basic properties of Gaussian processes.
A random function $\Xi$ is a centered Gaussian process if for any
fixed set $\{\theta_1, \ldots, \theta_n\}$ the random vector
$(\Xi(\theta_1), \ldots, \Xi(\theta_n))$ follows a multi-variate
normal distribution with mean zero.
The function $C(\theta, \theta') = \Cov(\Xi(\theta), \Xi(\theta'))$ is referred to as the covariance function.
Sampling a Gaussian process at points $\{\theta_1, \ldots, \theta_n\}$ is then as straightforward as drawing a Gaussian vector with mean zero and covariance matrix $\{C(\theta_i, \theta_j)\}_{i, j}$.

In this work, we make the common assumption that the covariance is stationary, i.e., it is only a
function of $\theta-\theta'$.
Thus, one can write $C(\theta, \theta') = k(\theta-\theta')$ for some function $k$.
The regularity of the random functions $\Xi$ is directly linked to the regularity of $k$; in this work we consider a classical squared exponential covariance function, 
\begin{equation}\label{eq:kernel}
    k(d) = \sigma^2 \exp\left(-\frac{d^2}{2l^2}\right),
\end{equation}
which results in rather smooth perturbations.
Here $\sigma>0$ controls the overall magnitude of the perturbations, and $l>0$ is a measure for its length scale.
For our problem, we cannot use the covariance function $k$ directly, because we need to guarantee periodicity of $\Xi$. To address this minor difficulty, we rely on the fact that any covariance function $k$ can be made periodic on $[0, 2\pi)$ by defining
\begin{equation}
    \tilde k(d) = \sum_{j \in \mathbb{Z}} k(d + j2\pi),
    \label{eq:periodiccovariance}
\end{equation}
see~\cite[(4.42)]{scholkopf_learning_2002}.
We thus use $\tilde{k}$ instead of $k$ for our optimization studies. For our construction of $\tilde{k}$, we truncate the sum (\ref{eq:periodiccovariance}) after just a few terms, since $k$ has exponential decay.

\begin{figure}[H]
   \begin{center}
       \includegraphics[height=4cm]{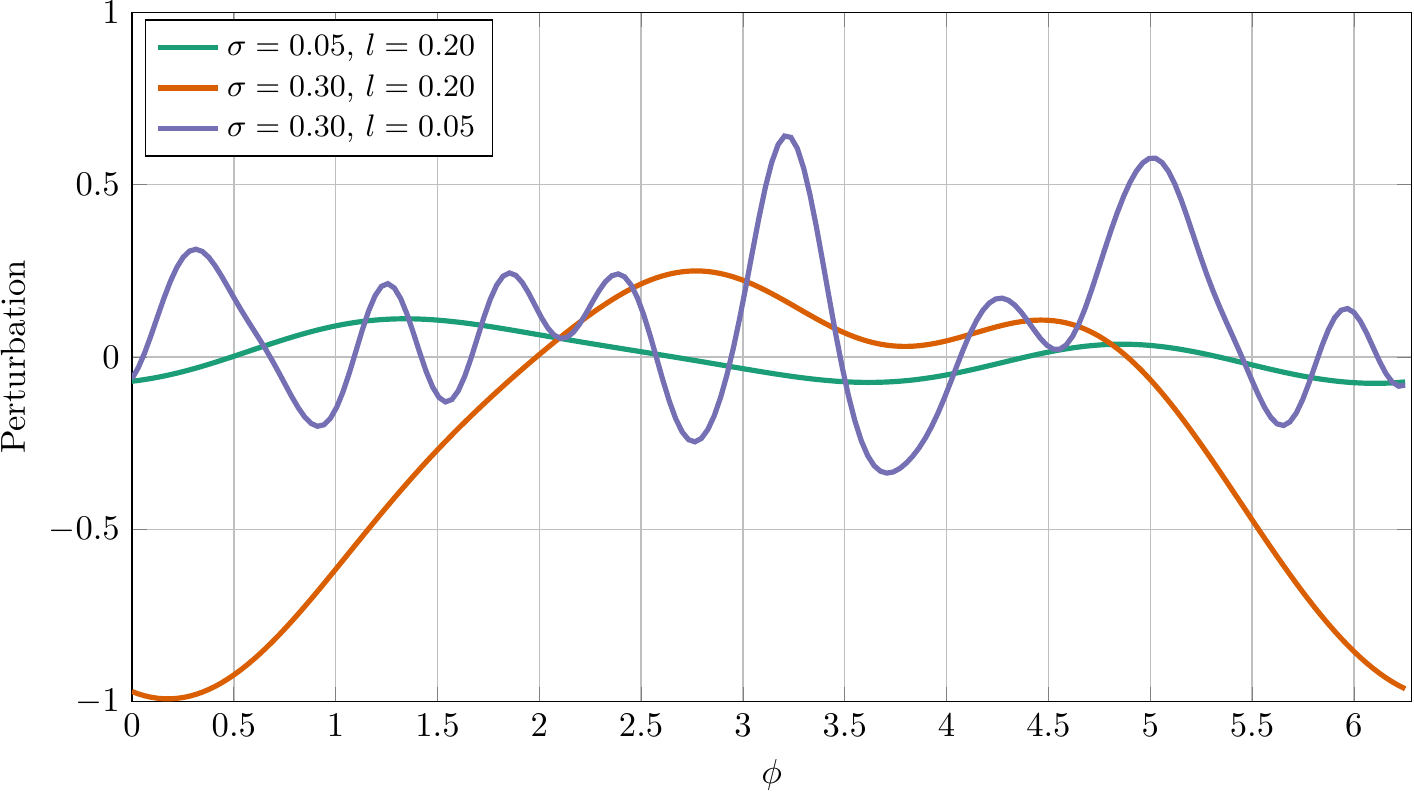}\hspace{1cm}
       \includegraphics[height=4cm]{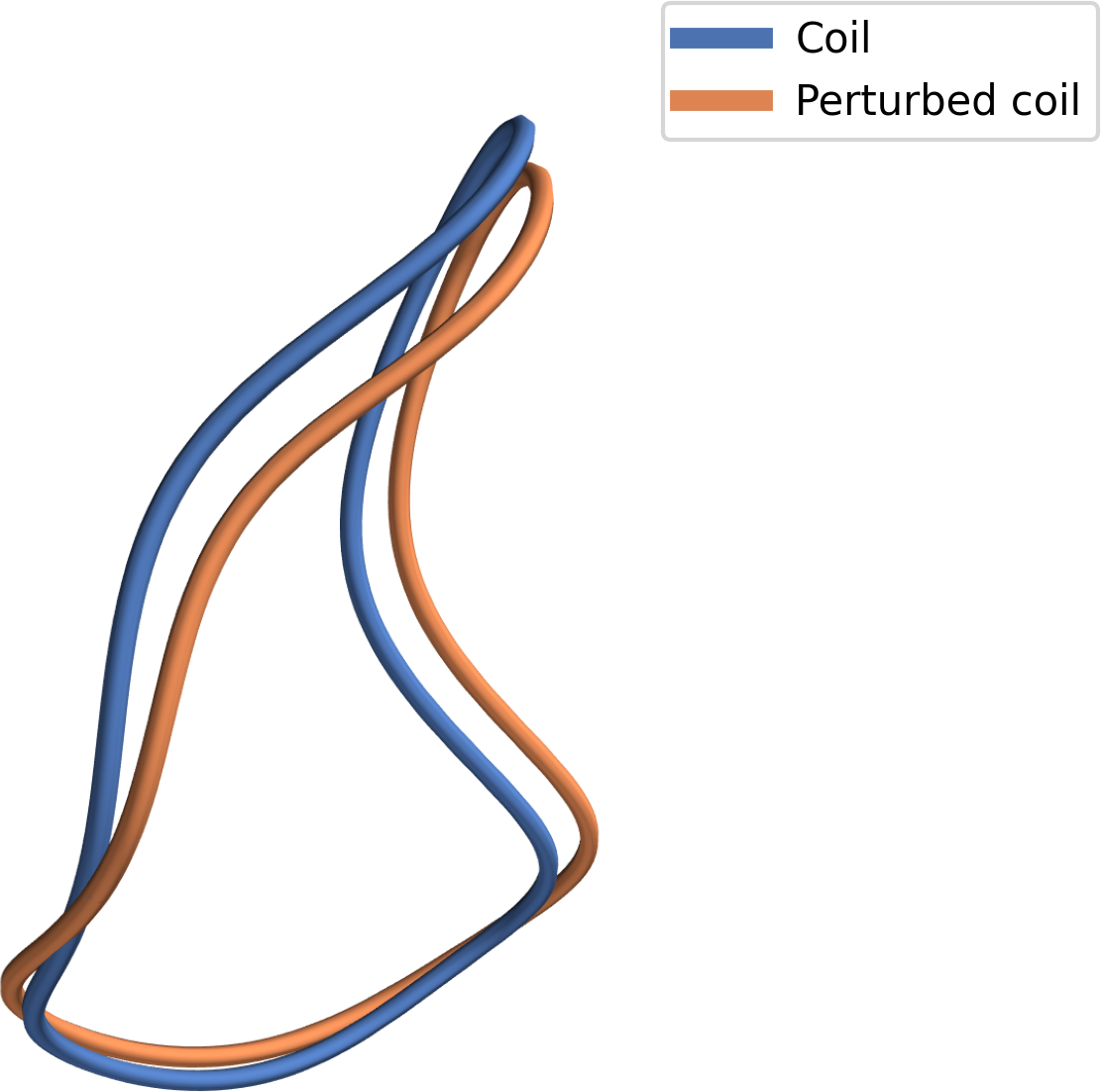}
   \end{center} 
\caption{Left: examples of periodic Gaussian processes for different
  parameter choices $\sigma$ and $l$ in \eqref{eq:kernel}. Right: a coil perturbed with a
  Gaussian process sample corresponding to $\sigma=0.1$,
  $l=0.2$. This large value of $\sigma$ has been chosen for
  illustration purposes. In our design process we choose smaller $\sigma$
  modeling realistic manufacturing errors.}\label{fig:gp-samples-and-pert}
\end{figure}

Finally, since our optimization algorithms require the knowledge of
derivatives, an additional property of Gaussian processes is relevant
to our work. By linearity
(see~\cite[\S9.4]{rasmussen_gaussian_2006},~\cite[\S2.2]{adler_geometry_2010})
the derivatives of $\Xi$ are also Gaussian and satisfy, for
$\theta,\bar\theta$:
\begin{equation}
    \begin{aligned}
        \Cov(\Xi'(\theta), \Xi(\bar\theta))  & = \partial_{\theta} C(\theta, \bar\theta) = \tilde{k}'(\theta-\bar\theta)\\
        \Cov(\Xi(\theta), \Xi'(\bar\theta))  & = \partial_{\bar\theta} C(\theta, \bar\theta) = -\tilde{k}'(\theta-\bar\theta)\\
        \Cov(\Xi'(\theta), \Xi'(\bar\theta)) & = \partial_{\theta}\partial_{\bar\theta} C(\theta, \bar\theta) = -\tilde{k}''(\theta-\bar\theta),
    \end{aligned}
\end{equation}
where the prime denotes derivative with respect to its argument. Thus, we can draw samples $\{\Xi(\theta_1), \ldots, \Xi(\theta_n),
\Xi'(\theta_1), \ldots, \Xi'(\theta_n)\}$ by drawing a Gaussian vector
with covariance matrix
\begin{equation}
    \begin{aligned}
    \Sigma
    &=\begin{bmatrix}
        \tilde{k}(\theta_1-\theta_1) & \ldots & \tilde{k}(\theta_1- \theta_n) & -\tilde{k}'(\theta_1-\theta_1) & \ldots &  -\tilde{k}^{'}(\theta_1-\theta_n) \\
        \vdots & \ddots & \vdots & \vdots & \ddots & \vdots \\ 
        \tilde{k}(\theta_n-\theta_1) & \ldots & \tilde{k}(\theta_n-\theta_n) & -\tilde{k}'(\theta_n-\theta_1) & \ldots & -\tilde{k}'(\theta_n-\theta_n) \\
        \tilde{k}'(\theta_1-\theta_1) & \ldots & \tilde{k}'(\theta_1-\theta_n) & -\tilde{k}''(\theta_1-\theta_1) & \ldots & -\tilde{k}''(\theta_1-\theta_n) \\
        \vdots & \ddots & \vdots & \vdots & \ddots & \vdots \\ 
        \tilde{k}'(\theta_n-\theta_1) & \ldots & \tilde{k}'(\theta_n-\theta_n) & -\tilde{k}''(\theta_n-\theta_1) & \ldots & -\tilde{k}''(\theta_n-\theta_n)
    \end{bmatrix}.
    \end{aligned}
\end{equation}
The standard approach to drawing such samples is to compute a matrix square root $\Sigma=LL^T$ (e.g.~via a Cholesky decomposition) and to draw a standard Gaussian vector $\zb$.
It is then straightforward to check that $\mathrm{Cov}(L\zb) = \Sigma$.
We show several random function draws from this Gaussian process as
well as a perturbed coil in Figure~\ref{fig:gp-samples-and-pert}.


\section{Stochastic and risk-averse optimization}\label{sec:stochastic}

In this Section, we briefly review the mathematical formulations of the different forms of stochastic optimization one may favour, depending on the coil optimization design goals and the level of risk one is willing to tolerate.

Let $f(\{\Gammab^{(i)}\}_{i=1}^{n_c}, \qb)$ be some quantity of interest which one wants to minimize, and which depends on the coil geometry $\{\Gammab^{(i)}\}_{i=1}^{n_c}$, as well as other quantities $\qb\in\re^{n_q}$ (e.g.~coil currents).
We define
\begin{equation}
    g((\{\Gammab^{(i)}\}_{i=1}^{n_c}, \qb), \{\Xib^{(i)}\}_{i=1}^{n_c}) = f(\{\Gammab^{(i)} + \Xib^{(i)}\}_{i=1}^{n_c}, \qb)
\end{equation}
For notational brevity we write $\xb = (\{\Gammab^{(i)}\}_{i=1}^{n_c}, \qb)$ and $\zetab = \{\Xib^{(i)}\}_{i=1}^{n_c}$, i.e., the variables to be optimized are contained in $\xb$ and the randomness is contained in $\zetab$.
For fixed $\xb$, $g(\xb, \zetab)$ is now a random variable, which needs to be scalarized in order to perform optimization.
Risk-neutral, risk-averse, and robust stochastic optimization formulations are all obtained by
different ways of scalarising $g(\xb, \zetab)$, and thus all take the distribution of manufacturing errors into account. While other stochastic
optimization formulations exist (e.g., \cite{ShapiroDentchevaRuszczynski09}),
we only describe these three below, since they are the most common formulations, and well suited to stellarator optimization. 

Before we do so, we observe that the deterministic optimization problem of minimizing $f$ is equivalent to minimizing
\begin{equation}\label{eqn:det-min}
    \min_{\xb} g(\xb, 0),
\end{equation}
i.e., it is assumed that no coil error is present. This is the approach traditionally taken in stellarator optimization, sometimes followed by perturbation tests at the optimal design 
\cite{Rummel_2004, zhu2017,Zhu_2018,giuliani2020} to evaluate the sensitivity of the objective with respect to coil errors. In contrast, the approaches
discussed next take the distribution
of manufacturing errors into account \textit{during} the optimization
process. 

\subsection{Types of optimization under uncertainty}
\subsubsection{Risk-neutral stochastic optimization}

Risk-neutral stochastic optimization corresponds to the situation in which one wants to find a solution $\xb$ that performs optimally with respect to the mean of the realisations. We thus obtain the optimization problem
\begin{equation}\label{eqn:riskneutral}
    \min_{\xb} \E{g(\xb, \zetab)},
\end{equation}
where $\E{\cdot}$ denotes expectation over the distribution of the
perturbations $\zetab$. In this article, we will mainly focus on this approach, and thus
sometimes simply refer to it as \emph{stochastic optimization}.
However, for certain stellarator design problems and certain objective functions, it can be desirable to explicitly avoid poor
objective values for some realizations of $\zetab$. This can be achieved
using risk-averse or robust formulations, which we summarize next.

\subsubsection{Risk-averse stochastic optimization: CVaR}
\begin{figure}[H]
   \begin{center}
       \includegraphics[width=0.99\textwidth]{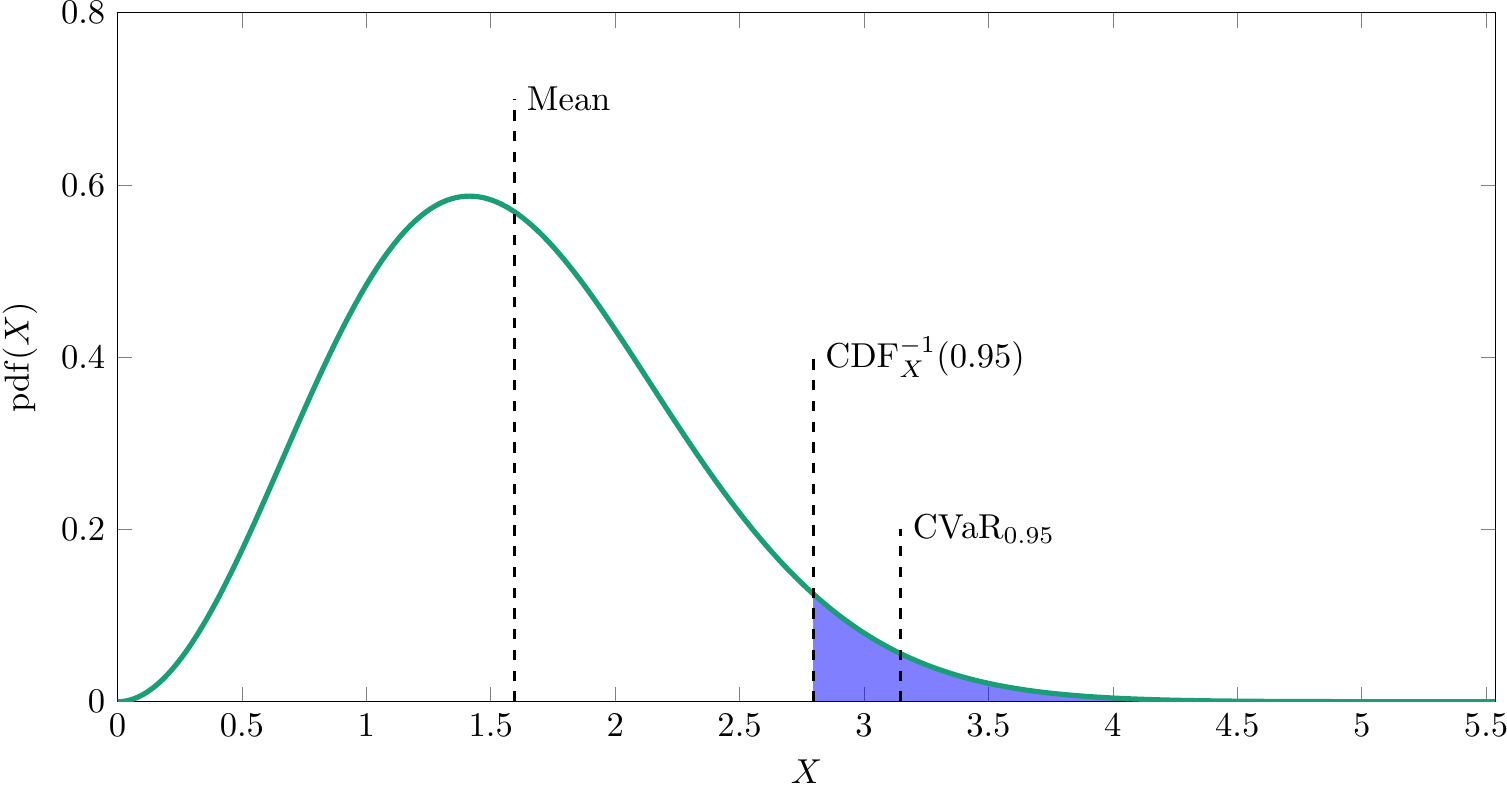}
   \end{center} 
   \caption{Illustration of mean and CVaR with $\alpha=0.95$ for
     an example distribution whose probability density function is
     shown in green.}\label{fig:distribution}
\end{figure}

A measure that focuses on the tail of the distribution is the
conditional value-at-risk, or CVaR.  The CVaR of a random variable $Z$
is defined as the expected value given that the random variable falls
into the $\alpha$-quantile of its distribution, i.e.,
\begin{equation}
    \cvar{\alpha}{Z} = \E{Z | Z>\mathrm{CDF}^{-1}_Z(\alpha)},
\end{equation}
where $\mathrm{CDF}_Z$ denotes the cumulative distribution function and
$\alpha\in[0,1]$. The difference between the risk-neutral formulation
and CVaR is illustrated with an example probability density function in 
Figure~\ref{fig:distribution}. This example highlights the fact that the
CVaR only depends on the tail of the distribution. One reason for the
popularity of CVaR over other risk-averse measures is its convexity as
well as the following equivalent formulation
(cf.~\cite[Theorem~1]{rockafellar_optimization_2000})
\begin{equation}
    \cvar{\alpha}{Z} = \inf_{t\in\re} \bigg[ t + \frac{1}{1-\alpha} \E{(Z-t)^+} \bigg],
\end{equation}
where $s^+:=\max(0,s)$, which is a convenient formulation for
computational purposes. To cope with the non-differentiable nature of the
$\max$-function, one relies in practice on a smooth approximation
$h_\eps:\re\to[0, \infty)$ with $h_\eps \to ({}\cdot{})^+$ as
$\eps\to 0$.
The specific form of $h_\eps$ used in this work is
\begin{equation}
    h_\eps(x) = \begin{cases}
        x & \text{ if } x \ge \eps/2,\\
        \frac{(x+\eps/2)^3}{\eps^2}-\frac{(x+\eps/2)^4}{2\eps^3} & \text{ if } -\eps/2 < x < \eps/2,\\
        0 & \text{ otherwise.}
    \end{cases}
\end{equation}
We thus obtain the optimization problem
\begin{equation}\label{eqn:cvarsmooth}
    \min_{\xb, t} \bigg[t + \frac{1}{1-\alpha} \E{h_\eps(g(\xb, \zetab)-t)}\bigg]
\end{equation}
for a small $\eps$. Risk-averse formulations based on CVaR are
successfully used in the insurance and finance industries and in
engineering
\cite{rockafellar_optimization_2000,KrokhmalPalmquistUryasev02,KouriSurowiec16} for instance.


\subsubsection{Robust stochastic optimization} 

To completely control the probability
of poor outcomes, one can optimize the worst possible scenario, i.e.,
\begin{equation}
    \min_{\xb} \max_{\zetab} g(\xb, \zetab).
\end{equation}
This formulation is typically combined with a model for randomness
that results in perturbations that are almost surely bounded. We will
observe in Section \ref{sec:num} that the difference between
risk-neutral and risk-averse coil designs is minor for the stellarator
optimization problem we consider. Thus, we do not explore robust
optimization further in this article, as it can be viewed as an extreme version of
risk-averse optimization.


\subsection{Sample average approximation}
The expected value in~\eqref{eqn:riskneutral}
and~\eqref{eqn:cvarsmooth} can typically not be computed analytically but has to be
approximated numerically.  For that purpose, we use the \emph{sample average approximation},
i.e., we
draw $N_\mathrm{MC}$ independent realisations of $\zetab_k$ and
approximate \begin{equation}\label{eq:SAA}
    \begin{aligned}
        \E{g(\xb, \zetab)} &\approx
        \frac{1}{N_\mathrm{MC}}\sum_{k=1}^{N_\mathrm{MC}} g(\xb,\zetab_k),\\
        \E{h_\eps(g(\xb, \zetab)-t)} &\approx
        \frac{1}{N_\mathrm{MC}}\sum_{k=1}^{N_\mathrm{MC}} h_\eps(g(\xb, \zetab_k)-t),
    \end{aligned}
\end{equation}
for the risk-neutral and risk-averse formulation respectively.
As $N_\mathrm{MC}\to \infty$, the random space approximation error is
of the typical Monte Carlo order
$O(N_\mathrm{MC}^{-1/2})$.  Since the samples $\zetab_k$ are kept fixed
throughout the optimization, \eqref{eq:SAA} results in a deterministic
optimization problem with $N_\mathrm{MC}$ terms. Note that by linearity the sample average
approximation of the gradients is exactly equal to the gradient of the sample
average approximation.

Our numerical tests in Section \ref{sec:num} focus on a comparison
between deterministic, risk-neutral and CVaR risk-averse stochastic
designs. Moreover, we study the role of the Monte Carlo sample size
$N_\mathrm{MC}$ for approximating the distribution.

\section{Direct coil design for quasi-symmetry in vacuum fields}\label{sec:single_stage}
In~\cite{giuliani2020} a new formulation was presented to directly design coils generating vacuum magnetic fields which are quasi-symmetric to high accuracy in a region close to the magnetic axis.
We briefly recall the basic structure of the objective that was developed there and then show the corresponding stochastic and risk-averse formulations.

Given a so-called expansion axis $\Gammab_\ab$ and real parameter $\bar\eta$ it was shown in~\cite{landreman_2018, landreman_2019} how to construct a magnetic field $\Bb_{QS}$ that is quasi-symmetric near the axis and how to compute its rotational transform $\iota$.
Calling $\Bb_{\text{coils}}$ the magnetic fields induced by the coils $\{\Gammab_\cb^{(i)}\}$, the approach of~\cite{giuliani2020} is then to find coils so that $\Bb_{QS} \approx \Bb_{\text{coils}}$.
Grouping the coefficients that describe the expansion axis and the real parameter $\eta$ in a vector $\ab$, and grouping the coefficients that describe the coils and their currents in a vector $\cb$, the objective is given by
\begin{equation}
    \begin{aligned}
        \hat J(\cb, \ab) &= \frac12 \int_{\Gammab_\ab} \|\Bb_\text{coils}(\cb)-\Bb_\text{QS}(\ab)\|^2 \dr l + \frac12 \int_{\Gammab_\ab} \|\nabla \Bb_\text{coils}(\cb)-\nabla \Bb_\text{QS}(\ab)\|^2 \dr l \\
                         &+ \frac12 \left( \frac{(\iota(\ab)-\iota_{0,a})^2}{\iota_{0,a}^2}\right) + R_a(\ab) + R_c(\cb),
    \end{aligned}
\end{equation}
where $\iota_{0,a}$ is a target rotational transform, and $R_a$ and $R_c$ contain various regularizations for the expansion axis and coils respectively.
The regularization terms include penalty functions for the length of the axis and the length of the coils, the curvature of the coils, and the distance between coils.
In~\cite{giuliani2020} it was shown that this formulation leads to an efficient method to design \emph{from scratch} coils producing nearly quasi-symmetric vacuum magnetic configurations, and to improve the quasi-symmetry properties of existing designs.

Random perturbations of the coils can be taken into account by
considering the objective
\begin{equation}
    \begin{aligned}
        \hat J(\cb, \ab, \zetab) &= \frac12 \int_{\Gammab_\ab} \|\Bb_\text{coils}(\cb, \zetab)-\Bb_\text{QS}(\ab)\|^2 \dr l + \frac12 \int_{\Gammab_\ab} \|\nabla \Bb_\text{coils}(\cb,\zetab)-\nabla \Bb_\text{QS}(\ab)\|^2 \dr l \\
                    &+ \frac12 \bigg( \frac{(\iota-\iota_{0,a})^2}{\iota_{0,a}^2}\bigg) + R_a(\ab) + R_c(\cb),
    \end{aligned}
\end{equation}
where $\Bb_\text{coils}(\cb, \zetab)$ corresponds to the magnetic field produced by the perturbed coils $\{\Gammab^{(i)} + \Xib^{(i)}\}$.
We emphasize here that the field $\Bb_\text{QS}$ and the rotational
transform are independent from the random variable $\zetab$, and hence we only
have to recompute $\Bb_\text{coils}(\cb, \zetab)$ and $\nabla \Bb_\text{coils}(\cb, \zetab)$ for different samples $\zetab$.
We observe that as it is stated here, the optimization problem can lead to numerical difficulties, because vanishing derivatives of $\Gammab^{(i)}$ result in a non differentiable curve length objective, and because of the large nullspace of the objective, which is partially due to the fact that different parametrizations give the same physical curve.
In order to address these numerical difficulties, we add the following regularization term, in addition to the axis length, coil length, coil curvature, and coil distance terms already included in~\cite{giuliani2020}, 
\begin{equation}
    R_{\mathrm{arc}}(\cb) = \sum_{i=1}^{n_c} \int_{[0,2\pi)} (\|{\Gammab^{(i)}}'(\theta)\|-t^{(i)})^2 d \theta,
\end{equation}
where $t^{(i)}=\frac{l^{i}}{2\pi}$ is the value that would correspond to a constant-arclength parametrization of a circle with coil length $l^{(i)}$.

\section{Numerical results for NCSX-like example}\label{sec:num}

\subsection{Implementation and setup}
We implement this optimization in the open source \pkg{PyPlasmaOpt} package available under \url{https://github.com/florianwechsung/PyPlasmaOpt}.
\pkg{PyPlasmaOpt} is a Python library that relies on the geometric objects and the Biot Savart implementation of the \pkg{SIMSOPT} stellarator optimization package \url{http://github.com/hiddenSymmetries/SIMSOPT}.
The implementation is parallelized across samples using MPI, and the Biot Savart computation is accelerated using SIMD instructions as well as OpenMP.
A function and gradient evaluation of a typical configuration (18 coils, 120 quadrature points per coil) with 1024 samples takes less than half of a second on a machine with two 24 Core Intel Xeon Platinum 8268 processors.

To solve the optimization problems
in~\eqref{eqn:det-min},~\eqref{eqn:riskneutral}, and~\eqref{eqn:cvarsmooth}, we
use the L-BFGS implementation in SciPy~\cite{virtanen2020scipy}.  For the
smoothed risk-averse objective, we use the risk-neutral minimizer as initial
guess, then solve the optimization problem, then reduce the smoothing parameter
$\eps$, and then solve the problem again, using the previous solution as
initial guess. This is repeated until $\eps=10^{-5}$.

We study a configuration that is inspired by the National Compact stellarator Experiment (NCSX).
NCSX consists of three distinct modular coils, which results in 18 coils after applying three fold rotational symmetry and stellarator symmetry.
We re-emphasize  that while the design space only consists of coils that satisfy these symmetries, the perturbations are not required to satisfy any form of symmetry.
For full details of the configuration we refer to \cite[Section~6]{giuliani2020}.
Finally, we note that NCSX was designed for finite pressure equilibria, as opposed to the vacuum field that we consider here, and that the additional planar coils of the NCSX design are ignored in our calculations.

We study this problem in particular as the NCSX project was cancelled
due to increasing costs because of, among other reasons, the
requirement of tight engineering tolerances on the coils.  To model
the distribution of manufacturing errors, we choose a length scale of
$l=0.4\pi$ and a standard deviations of either $\sigma = 10^{-2}$ or
$\sigma=3\times 10^{-3}$ in the kernel \eqref{eq:kernel} defining the
Gaussian process. For each of the coils, these values correspond to manufacturing errors of a few centimeters for $\sigma=10^{-2}$ and several millimeters for $\sigma=3\times 10^{-3}$ .

\begin{figure}
    \begin{center}
        \includegraphics[trim={4cm 2.5cm 2cm 2.3cm},clip, width=0.49\textwidth]{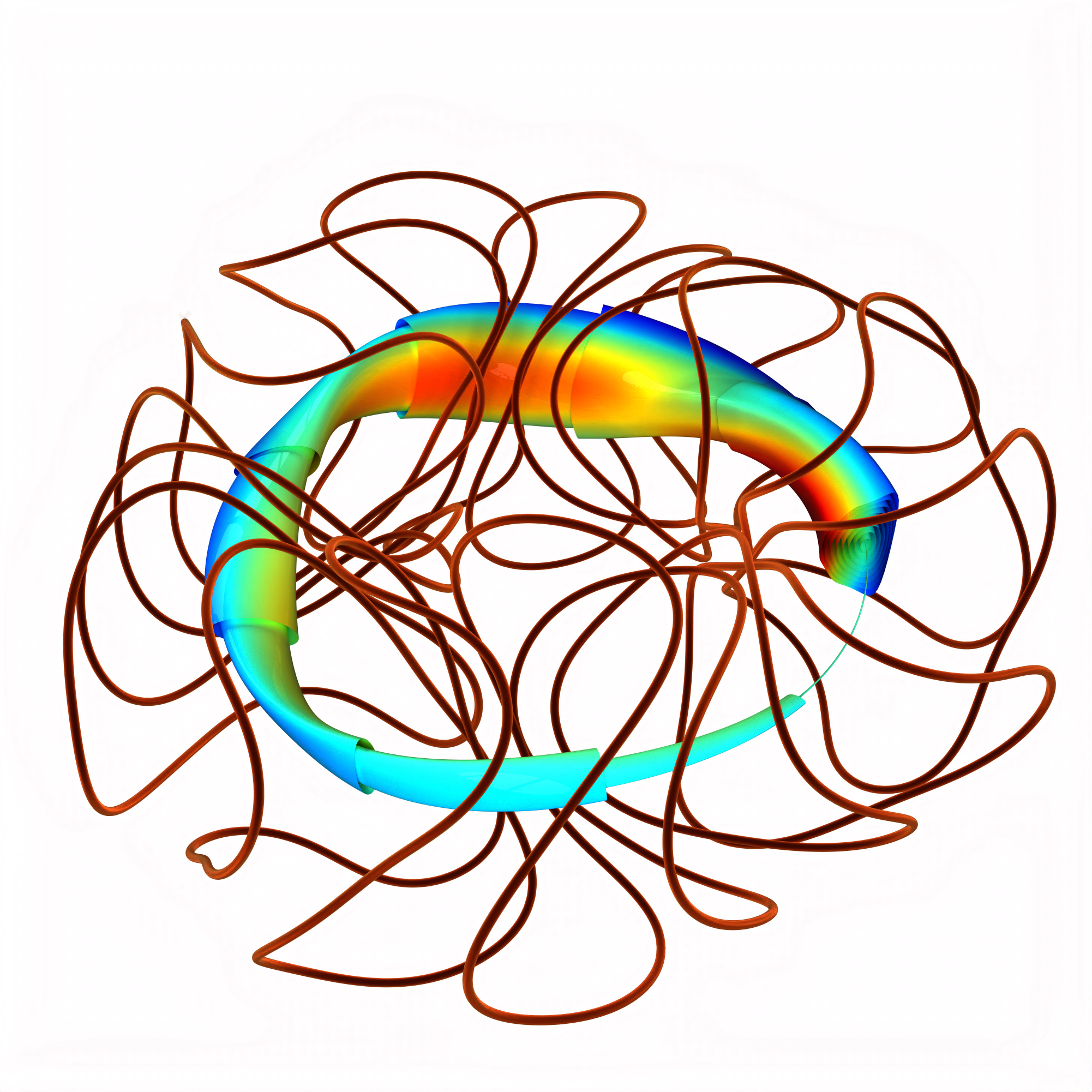}
        \includegraphics[trim={4cm 2.5cm 2cm 2.3cm},clip, width=0.49\textwidth]{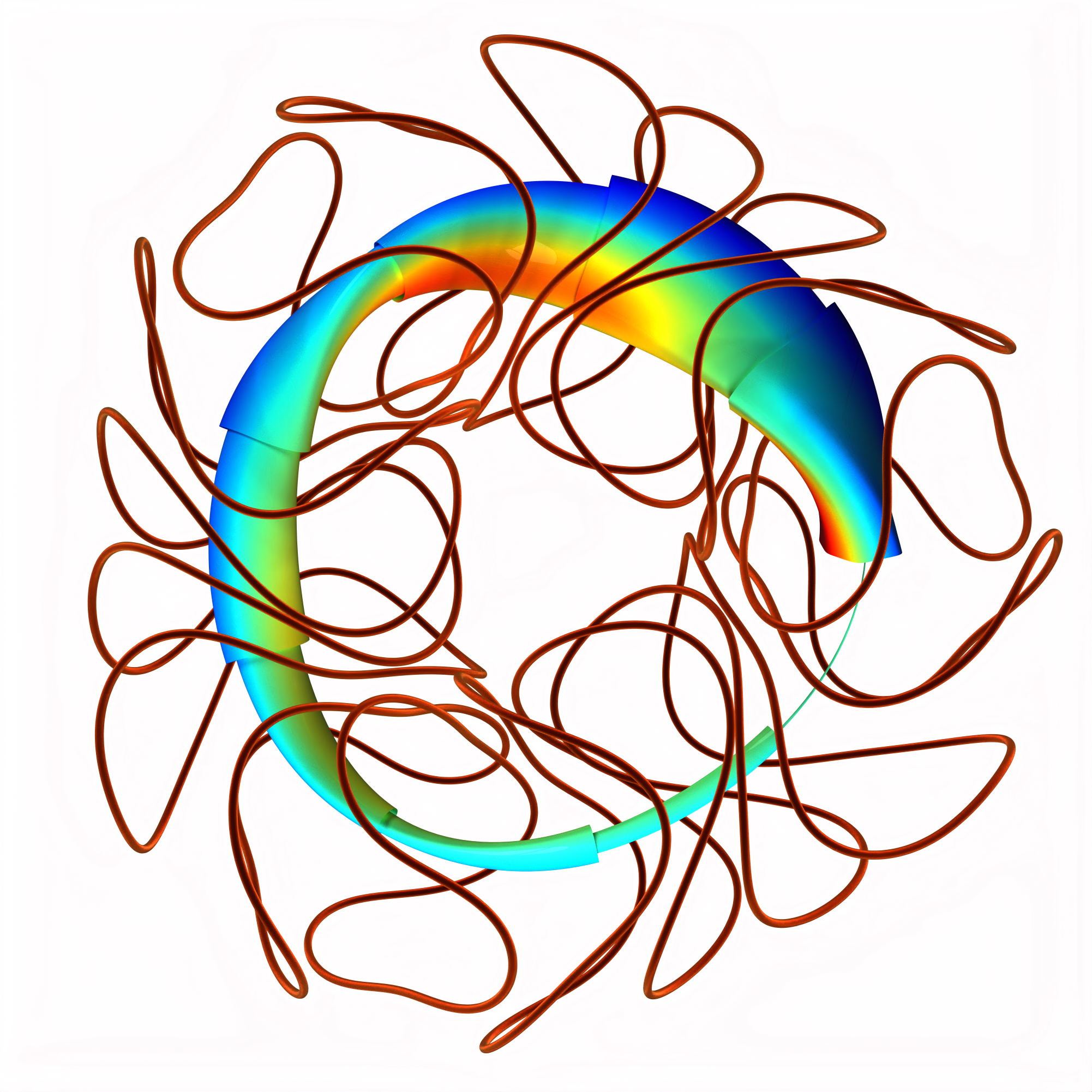}
    \end{center}
    \caption{Two views of an optimal stellarator design computed using
      risk-neutral stochastic optimization with $N_{\mathrm{MC}}=1024$ samples with
      $\sigma=10^{-2}$.  We also show a range of nested magnetic
      surfaces, with cool colors corresponding to low field strength
      and warm colors corresponding to high field
      strength.}\label{fig:ncsx-opt}
\end{figure}

When running the optimization algorithm with different initial
guesses, we observe that multiple local minima exist.  For this
reason, in the following sections we usually show results for
several different minima, that were obtained by starting the optimizer
from eight different initial guesses.  These initial guesses were
obtained by randomly perturbing the Fourier coefficients describing
the coils using independent normal random variables with standard
deviation of $0.01$. 
One of the obtained configurations together with a range of magnetic
surfaces is shown in Figure~\ref{fig:ncsx-opt}.
Next, in Sections \ref{subsec:determ} and \ref{subsec:out-of-sample}
we mainly focus on the role the number of samples $N_\textrm{MC}$ in
the sample average approximation \eqref{eq:SAA} plays in approximating
the distribution of coil perturbations, and on the sensitivity of both the deterministic and the stochastic minimizers on the initial guesses for the optimization algorithm.

\subsection{Deterministic versus stochastic designs}\label{subsec:determ}
\begin{figure}[H]
    \begin{center}
        \includegraphics[width=0.32\textwidth]{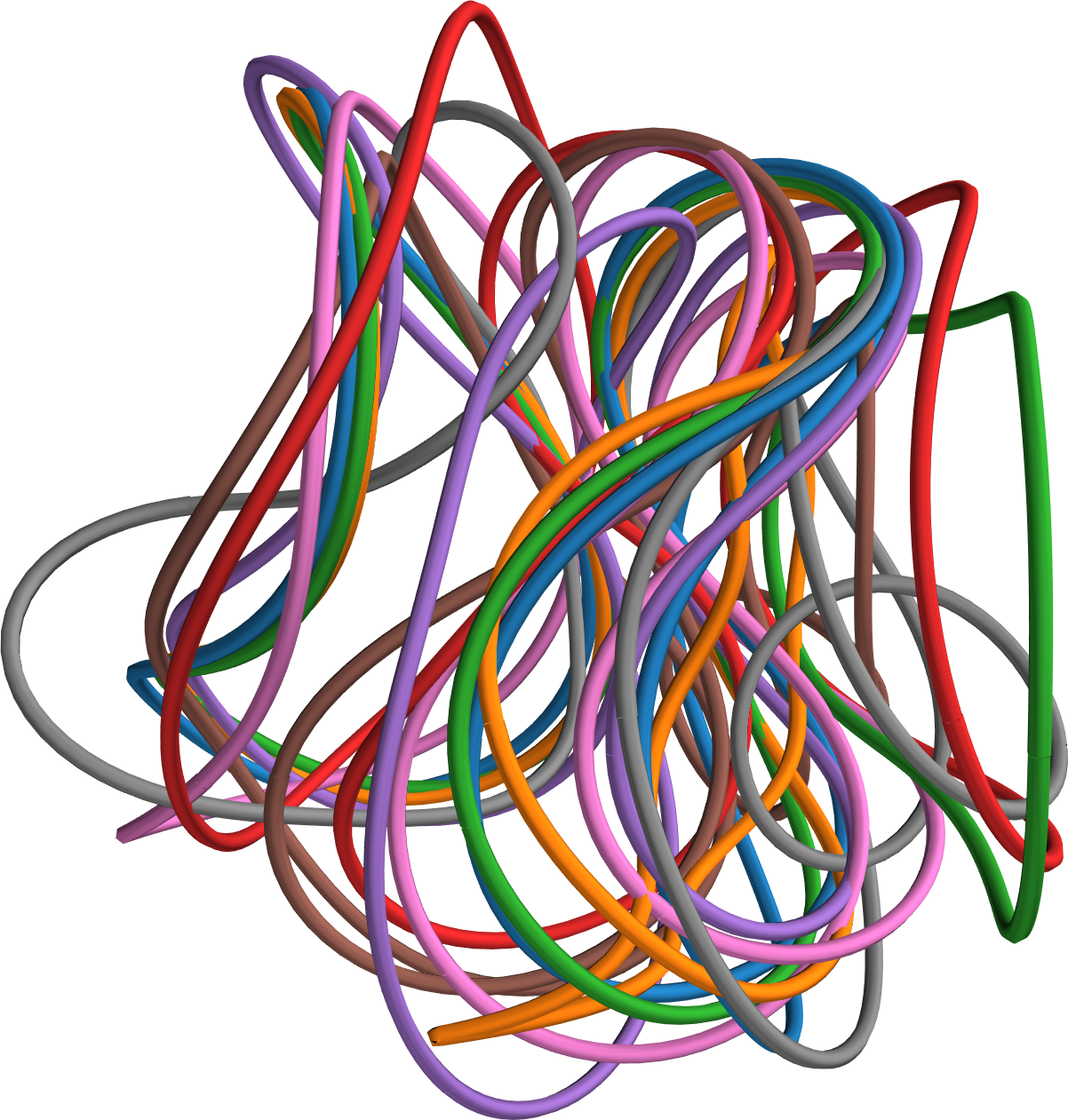}
        \includegraphics[width=0.32\textwidth]{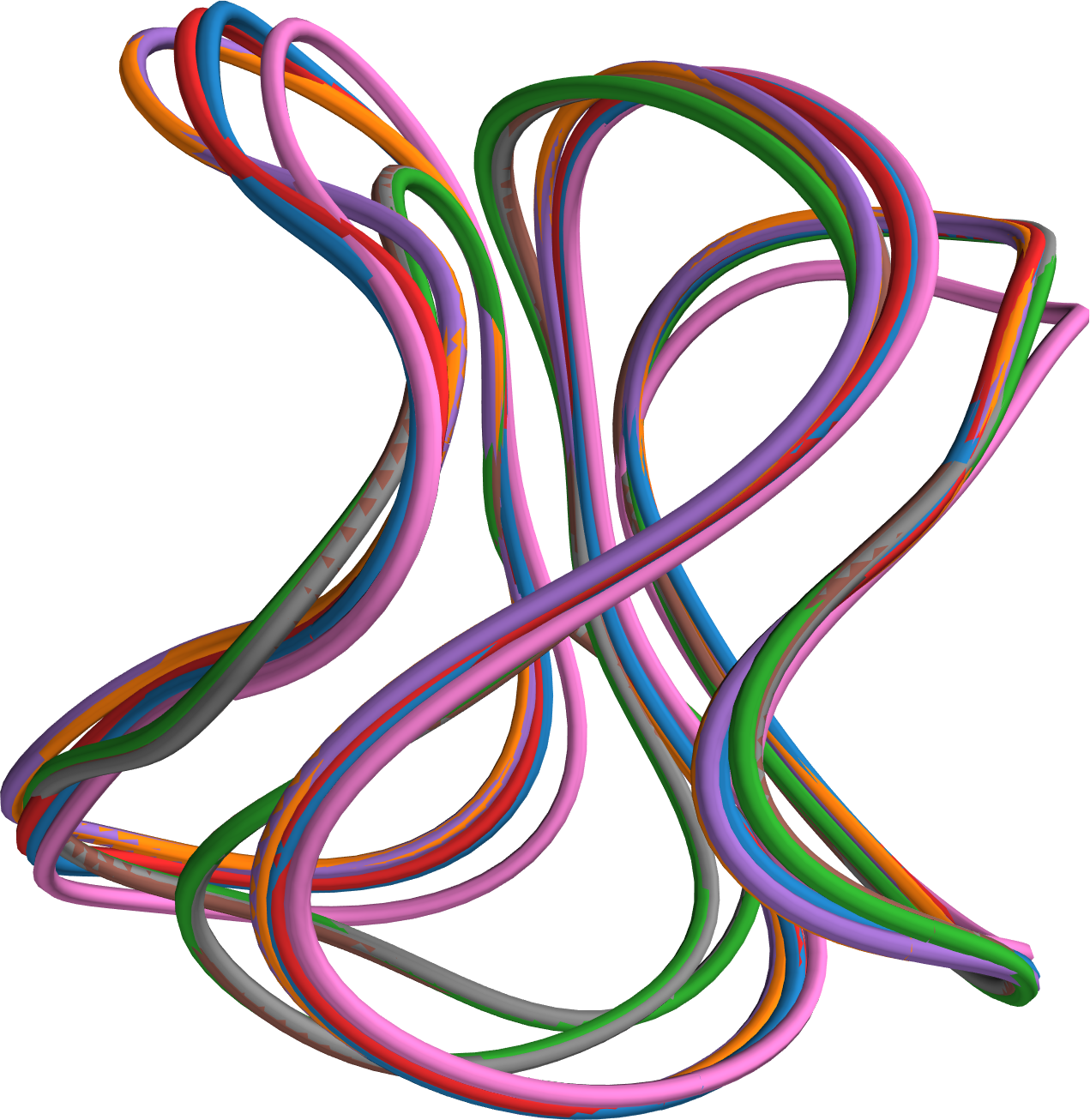}
        \includegraphics[width=0.32\textwidth]{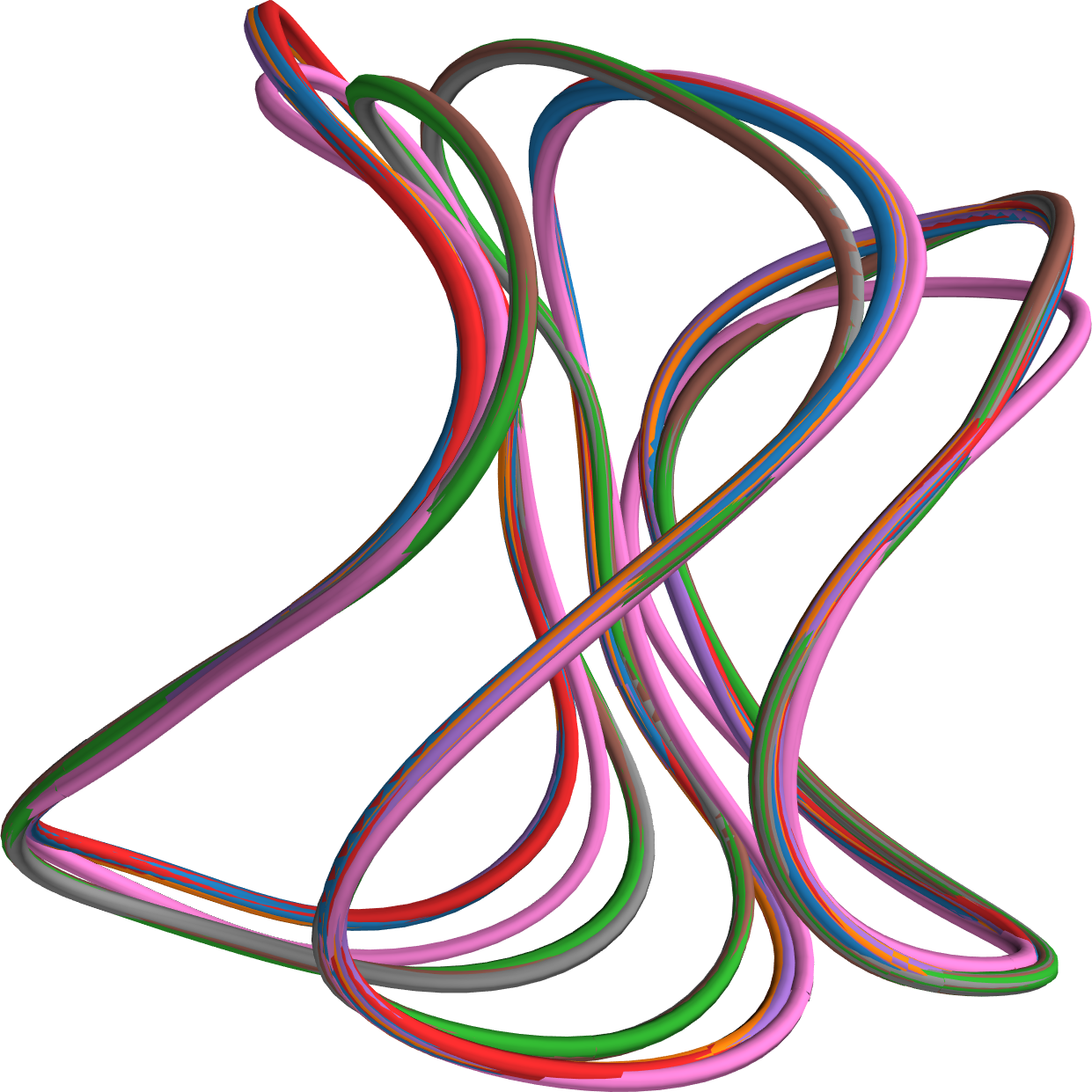}
    \end{center}
    \caption{Optimal coil designs for the three independent coils. The
      designs are obtained by using deterministic optimization (left),
      risk-neutral stochastic optimization with $N_{\mathrm{MC}}=4$ samples (middle) and
      with $N_{\mathrm{MC}}=1024$ samples (right) for eight different initial guesses
      (i.e.,~each panel shows 24 coil shapes).  Different colors
      correspond to different minimizers.}\label{fig:ncsx-minimizers}
\end{figure}

Figure~\ref{fig:ncsx-minimizers} shows minimizers that were obtained
by solving the deterministic and the stochastic optimization problem
for eight different initial guesses.  We clearly see that for the
deterministic problem there exists a large number of different
minimizers that are all distinct from each other.  As coil errors
(with $\sigma=10^{-2}$) are introduced and the sample average
    approximation size $N_\textrm{MC}$ (see \eqref{eq:SAA}) is
    increased, the number of distinct minimizers that the algorithm
    finds is reduced and there is less variation between the different coil designs.
    We note that the gradient was reduced by over 10 orders of
    magnitude for all of these minimizers, so we can be confident that this result is not an artifact of a possibly incomplete convergence of the algorithm. It is remarkable that stochastic design formulations
    are significantly less prone to having multiple
    minima. Intuitively, this can be understood by the fact that
    stochastic designs must perform well on average for a distribution
    of manufacturing errors, which prevents overfitting and also makes
    the objective locally convex in most directions around a
    minimizer. In the following discussion, we find several additional
    advantages that designs computed from a stochastic formulation
    have compared to designs based on a deterministic formulation.

\subsection{Out of sample distribution at the minimizer}\label{subsec:out-of-sample}
\begin{figure}[H]
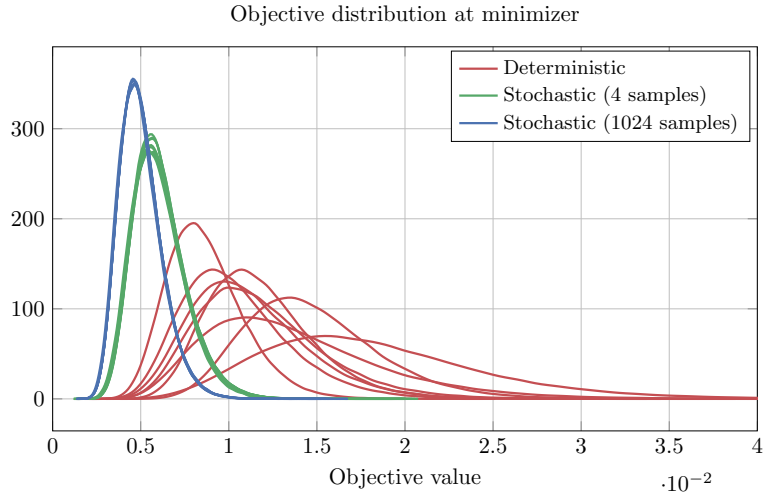

    \includestandalone[width=\textwidth]{./plots/density}
    \caption{Kernel density estimate for the distribution of the objective value evaluated at
      minimizers obtained by running deterministic (red) and
      stochastic optimization with $N_{\mathrm{MC}}=4$ and $N_{\mathrm{MC}}=1024$ samples used in the sample average approximation (green and blue,
      respectively).  Eight different initial guesses are used for
      each optimization, leading to different minimizers and thus
      different distributions. The distribution is approximated using
      $2^{18}=262,144$ independent samples drawn from the
      error distribution.}\label{fig:ncsx-distribution}
\end{figure}

As we only optimized the mean of the objective for a finite number of perturbations, we have to check that this performance generalises to the full distribution of coil errors.
To do this we draw a large number ($2^{18}$) of new samples, which are different from those used in the sample average approximation, and evaluate the objective at the minimizer for each of the samples.
In the context of statistics and machine learning, this procedure is known as out-of-sample testing or cross validation.
Figure~\ref{fig:ncsx-distribution} shows the resulting distribution for minimizers obtained from deterministic and stochastic optimization from eight different initial guesses.
We see that the performance of the minimizers of the deterministic
optimization problem varies strongly on perturbed coils. 
The minimizers found by stochastic optimization have lower objective
value on average and for most perturbations. In other words, the stochastic designs perform significantly
better than the deterministic designs.
Additionally, different minimizers obtained as a
result of different initializations of the algorithm perform vastly
differently for the deterministic formulation, but very similarly for the
stochastic formulations, in particular when 1024 samples are used to estimate the expected value in~\eqref{eq:SAA} for the stochastic formulation.

\subsection{Quasi-symmetry close to the axis}\label{subsec:QS}
\begin{figure}[H]
    \begin{center}
    \includegraphics[width=0.32\textwidth]{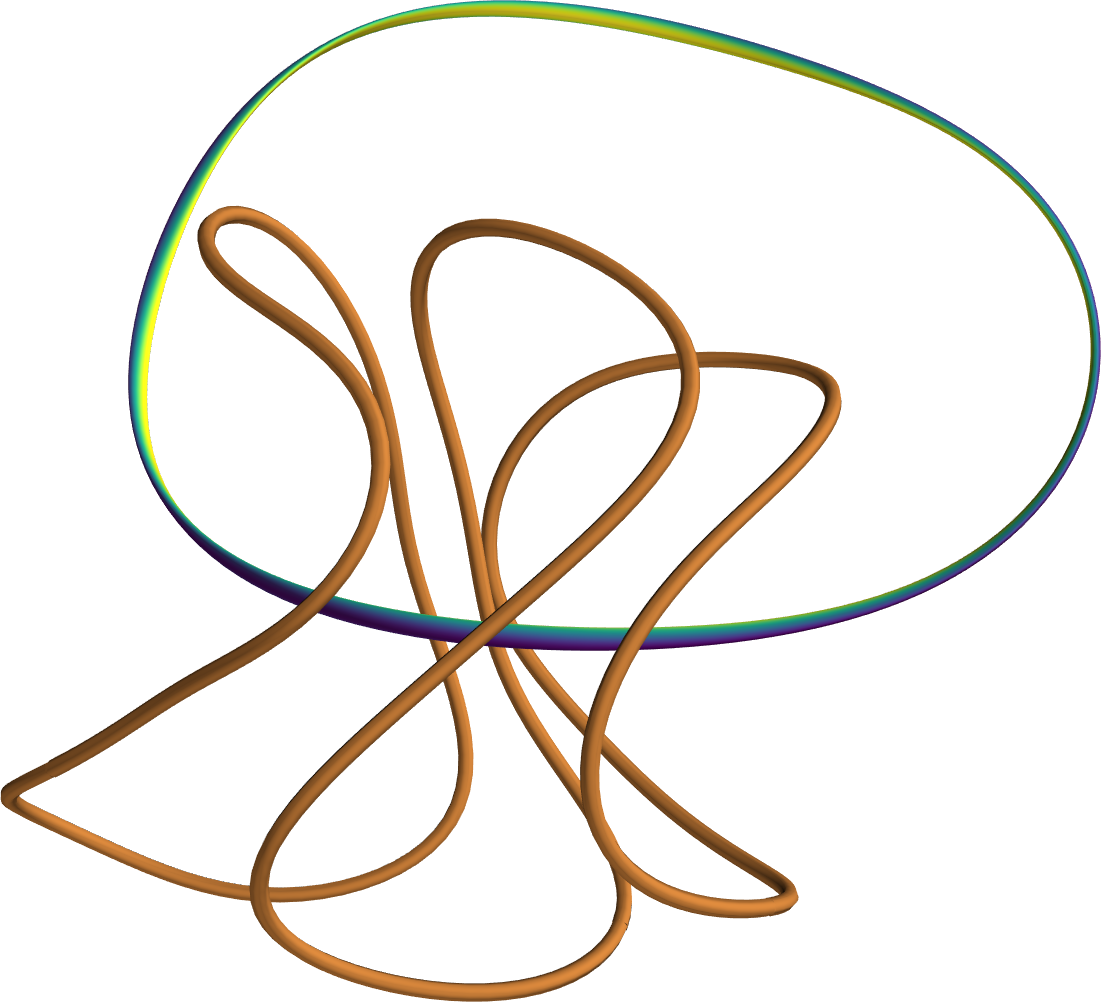}
    \includegraphics[width=0.32\textwidth]{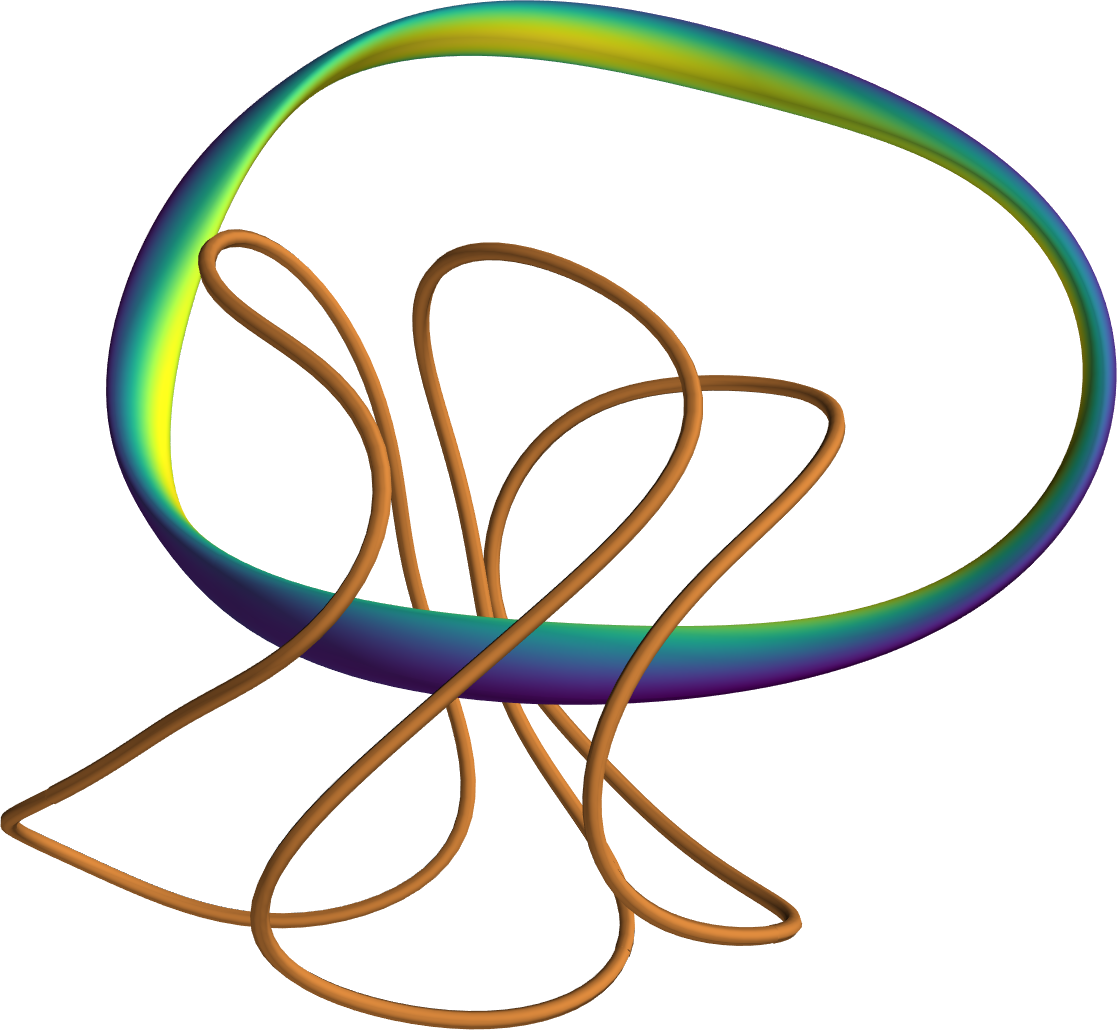}
    \includegraphics[width=0.32\textwidth]{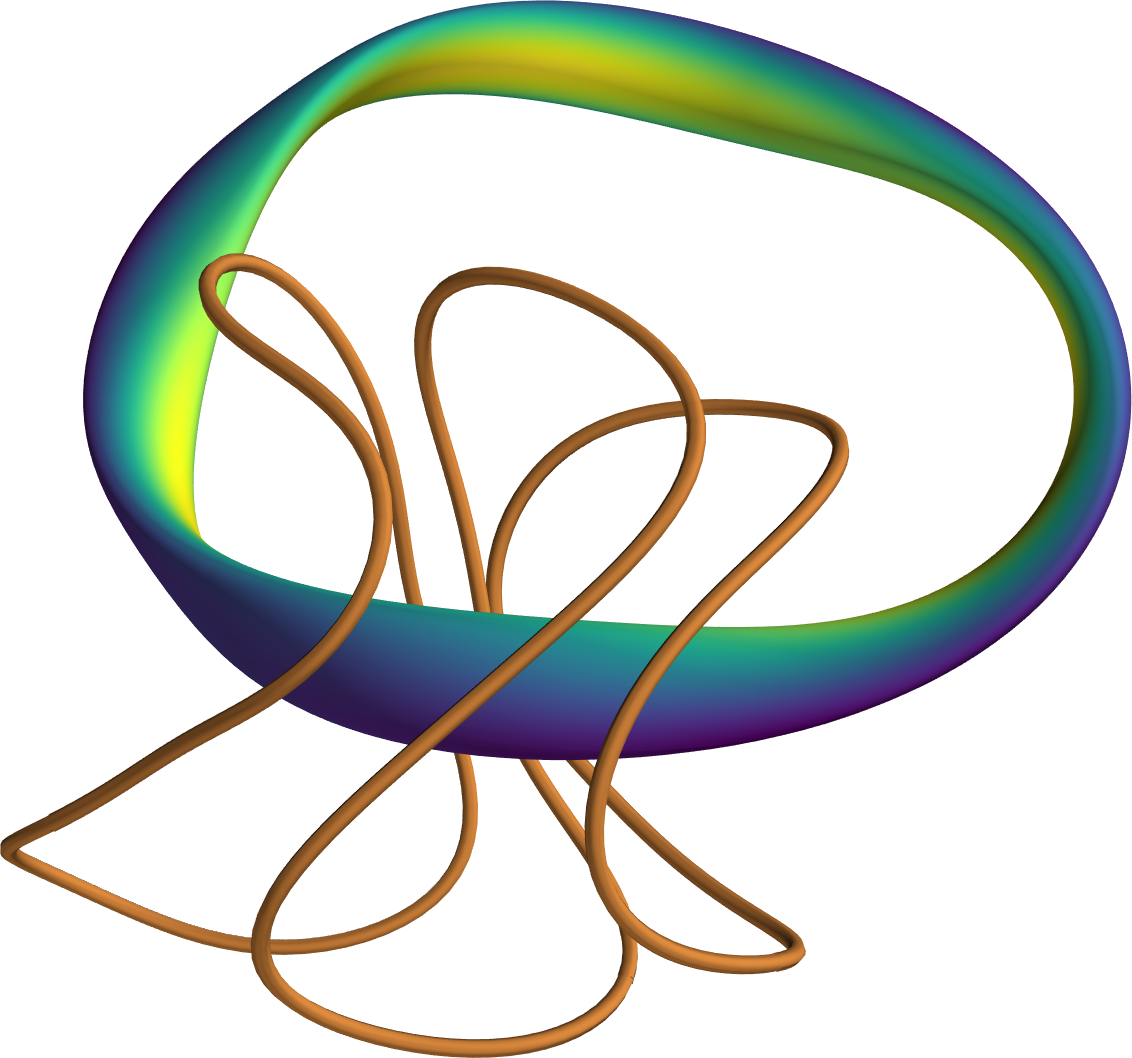}\vspace{0.5cm}
    \includegraphics[width=0.32\textwidth]{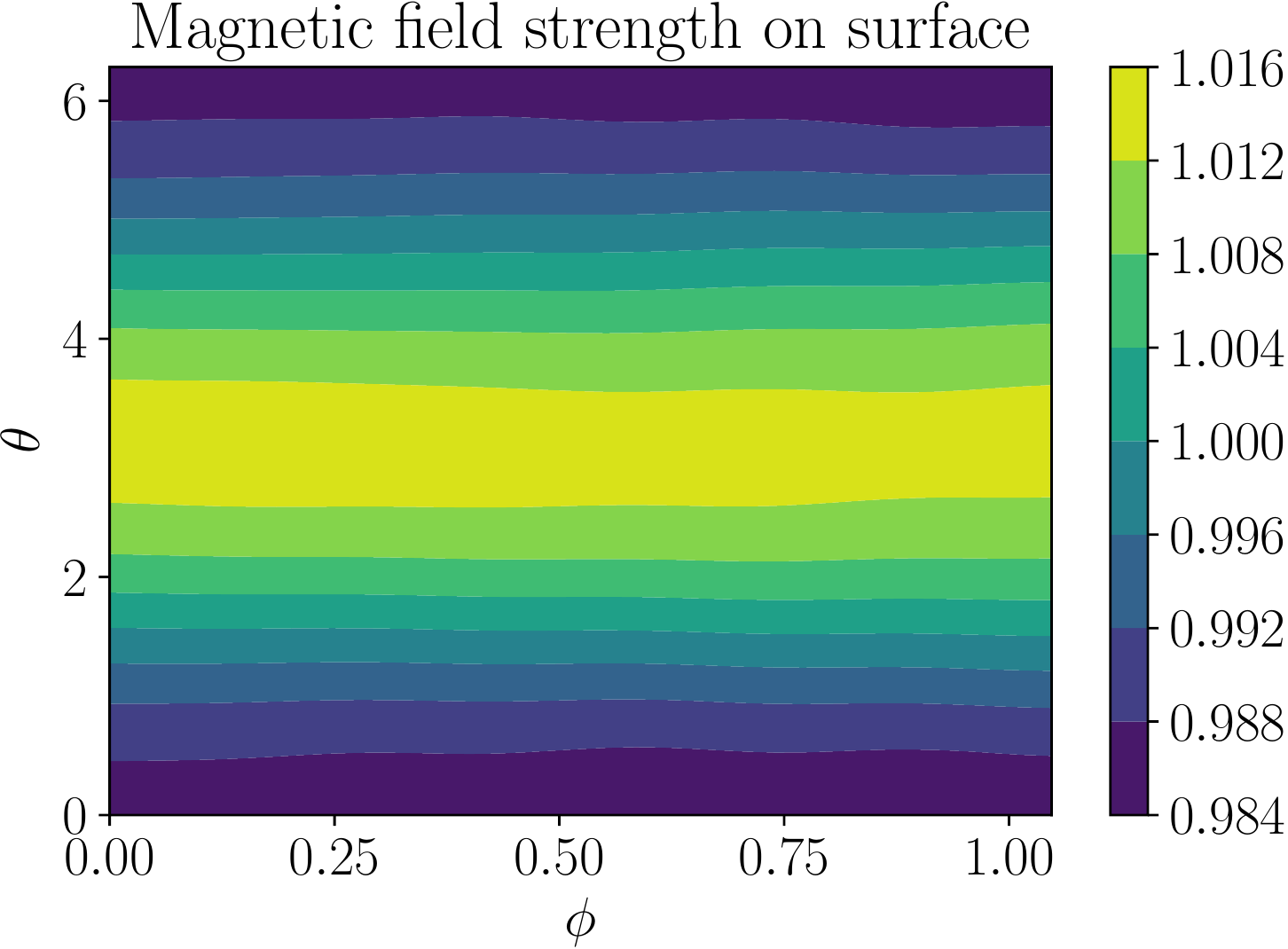}
    \includegraphics[width=0.32\textwidth]{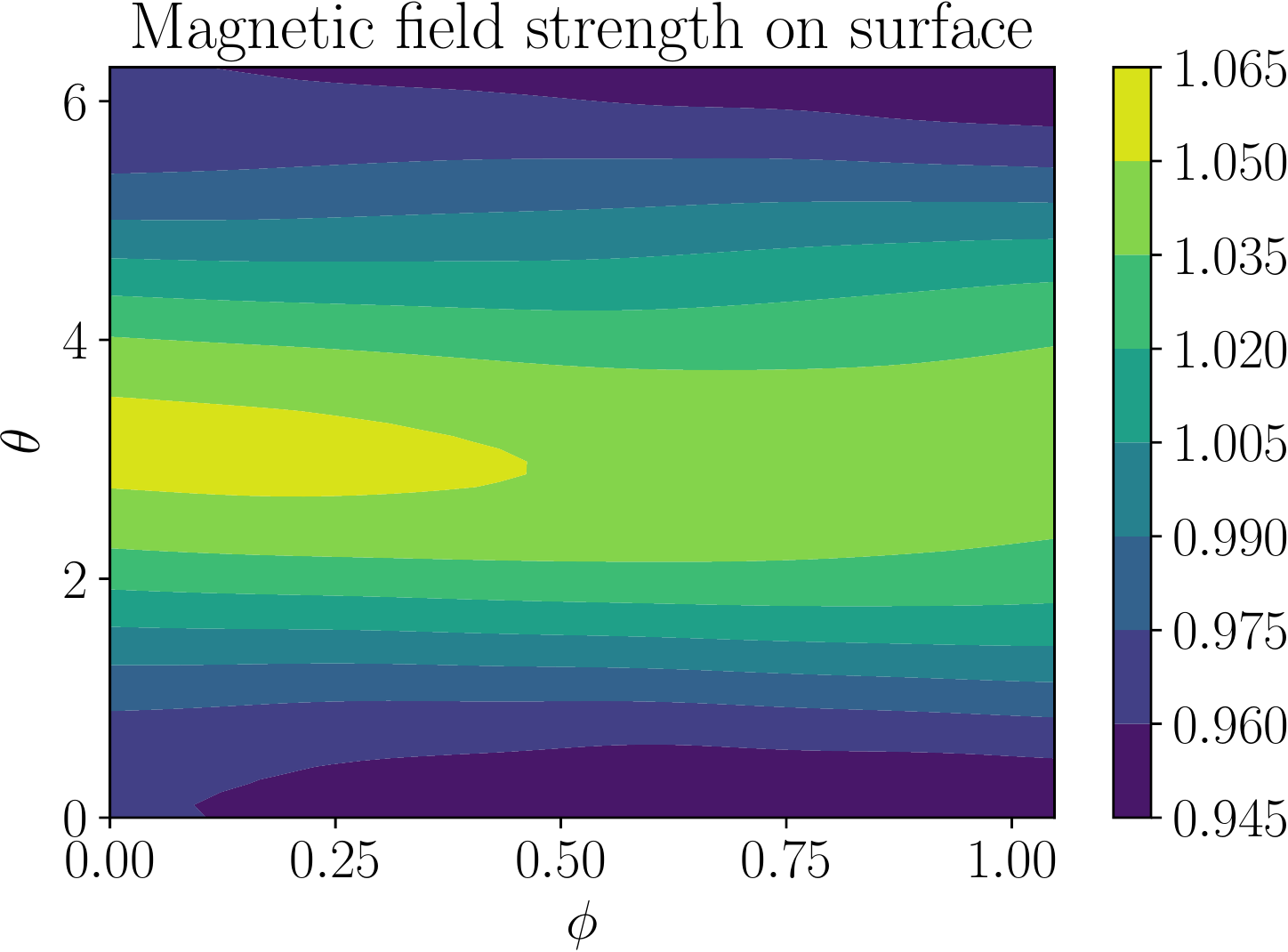}
    \includegraphics[width=0.32\textwidth]{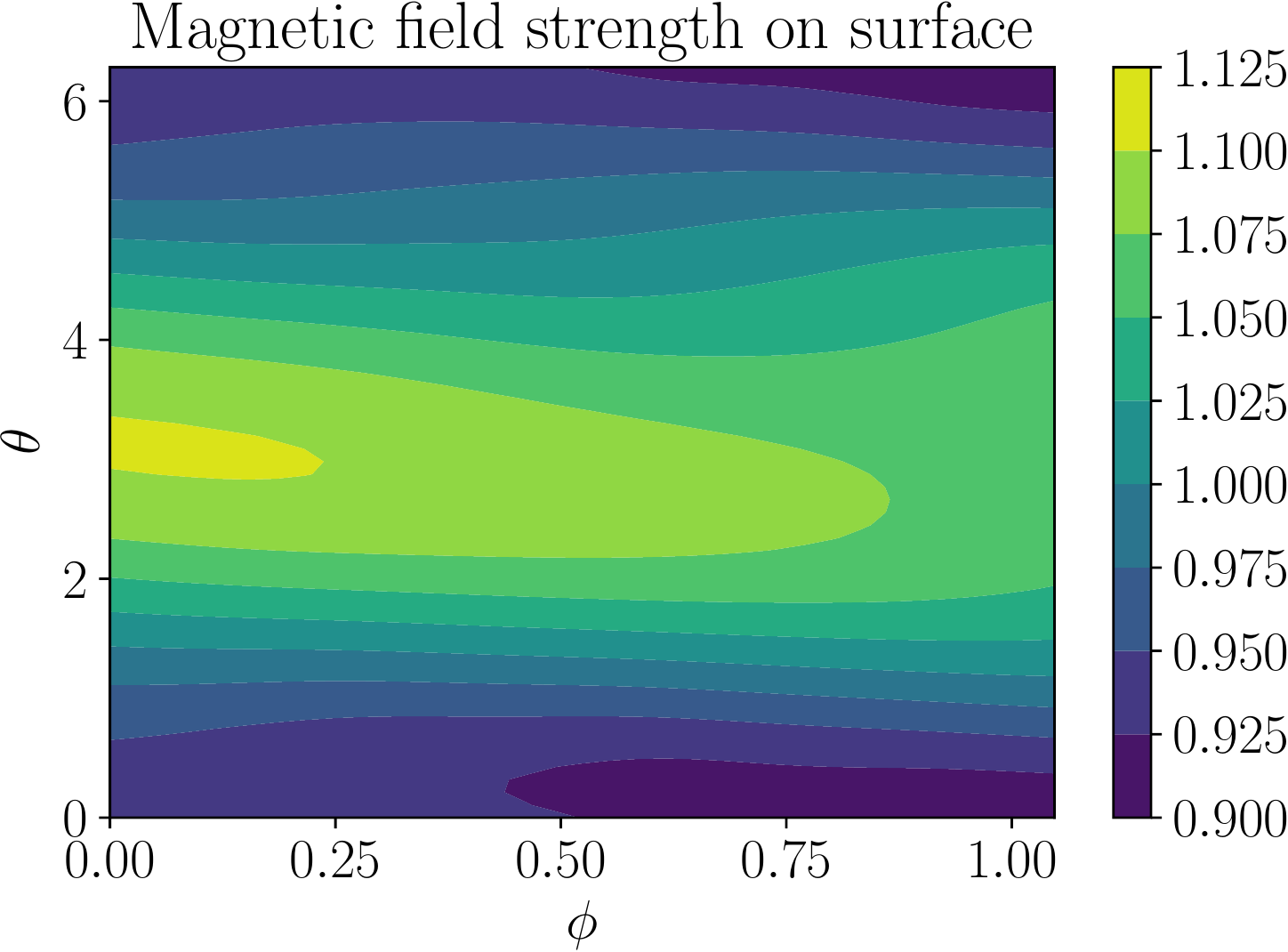}
    \caption{Top: Coils and magnetic surfaces in Boozer coordinates for the best
      configuration obtained from stochastic optimization with $N_{\mathrm{MC}}=1024$
  samples. Bottom: Strength of the magnetic field as a function of
  Boozer angles $\phi$ and $\theta$ on the surfaces shown on the top.
  Perfect quasi-symmetry corresponds to the magnetic field strength being independent of $\phi$.}\label{fig:ncsx-surfaces}
    \end{center}
\end{figure}

The objective is designed to ensure quasi-symmetry near the axis.
To confirm that this is achieved and to investigate the magnetic field away from the axis, we compute magnetic flux surfaces $\Sb(\phi,\theta)$ parametrized by Boozer angles $\phi$ and $\theta$.
We recall that a magnetic field is called quasi-axisymmetric if $|\Bb(\Sb(\phi, \theta))|$ is a function of $\theta$ only.
Figure~\ref{fig:ncsx-surfaces} shows three surfaces computed for the best configuration obtained from stochastic optimization with $N_{\mathrm{MC}}=1024$ samples.
We can see that for the surface closest to the axis, the field strength is indeed close to constant in $\phi$.
As we move away from the axis, this property is gradually lost.

\begin{figure}[H]
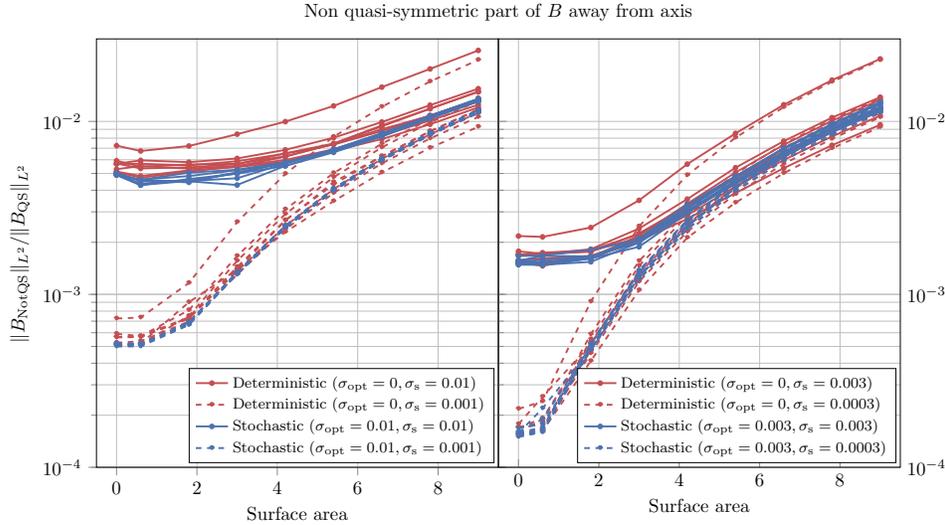

    \includestandalone[width=\textwidth]{./plots/qs_on_surface}
    \caption{Mean non-quasi-symmetry on a range of surfaces for eight different minimizers obtained from deterministic and stochastic ($N_{\mathrm{MC}}=1024$ samples) optimization. Here we use the surface area as a label for the surfaces. Larger surface areas correspond to surfaces farther from the axis.
    }\label{fig:ncsx-non-qs}
\end{figure}

To quantify this statement, we decompose $|\Bb|$ into a quasi-symmetric and non-quasi-symmetric part by defining
\begin{equation*}
    \begin{aligned}
        \|\Bb\|_{\mathrm{QS}}(\theta) &= \frac{\int_{[0,2\pi)} \|\Bb(\Sb(\phi, \theta))\| \|\partial_\phi \Sb \times \partial_\theta\Sb\|d\phi}{\int_{[0,2\pi)} \|\partial_\phi \Sb \times \partial_\theta\Sb\|d\phi}\\
        \|\Bb\|_{\mathrm{NotQS}}(\phi, \theta) &= \|\Bb(\Sb(\phi, \theta))\| - \|\Bb\|_{\mathrm{QS}}(\theta).
    \end{aligned}
\end{equation*}
We then measure the norm of the non-quasi-symmetric part and report $\big[\frac{\iint_{\Sb} \|\Bb\|_{\mathrm{NotQS}}^2 dS}{\iint_{\Sb} \|\Bb\|_{\mathrm{QS}}^2 dS}\big]^{1/2}$ in Figure~\ref{fig:ncsx-non-qs}.
Since we only enforce quasi-symmetry near the axis, we expect this
measure to be small for surfaces close to the axis, and to increase as we move away from the axis.

We perform the stochastic optimization for $\sigma_\mathrm{opt}=0.01$ and $\sigma_\mathrm{opt}=0.003$.
Quasi-symmetry is then evaluated by drawing 20 new samples with standard deviation $\sigma_s$ and computing surfaces for the fields induced by the perturbed coils.
The case of $\sigma_s < \sigma_\mathrm{opt}$ can be viewed as the estimate for the perturbation size in the original optimization being pessimistic, or be due to improvements in the manufacturing process between design and construction of the coils.
For comparison, we also compute surface and non-quasi-symmetry for the minimizers obtained from deterministic optimization (corresponding to $\sigma_\mathrm{opt}=0$).
The results are displayed in Figure~\ref{fig:ncsx-non-qs}.

Overall we observe that the configurations have very little non-quasi-symmetric contribution close to the axis, and that the non-quasi-symmetry then increases away from axis.
As can be expected, the difference between stochastic and deterministic optimization is most significant for large coil perturbations (left plot, solid lines).
Keeping those same configurations, but reducing the perturbation magnitude in the newly drawn samples, we can see that the configurations obtained from both types of optimization benefit from the added accuracy (left plot, dashed lines).
Importantly, the designs from stochastic optimization systematically perform better than those from deterministic optimization.

For smaller perturbations (right) the behaviour remains qualitatively the same, but the overall difference between stochastic and deterministic optimization is smaller.
This is expected, as in the limit of $\sigma_\mathrm{opt}\to0$ the
two optimization strategies become identical.

\subsection{Rotational transform on axis}\label{subsec:iota}
\begin{figure}[H]
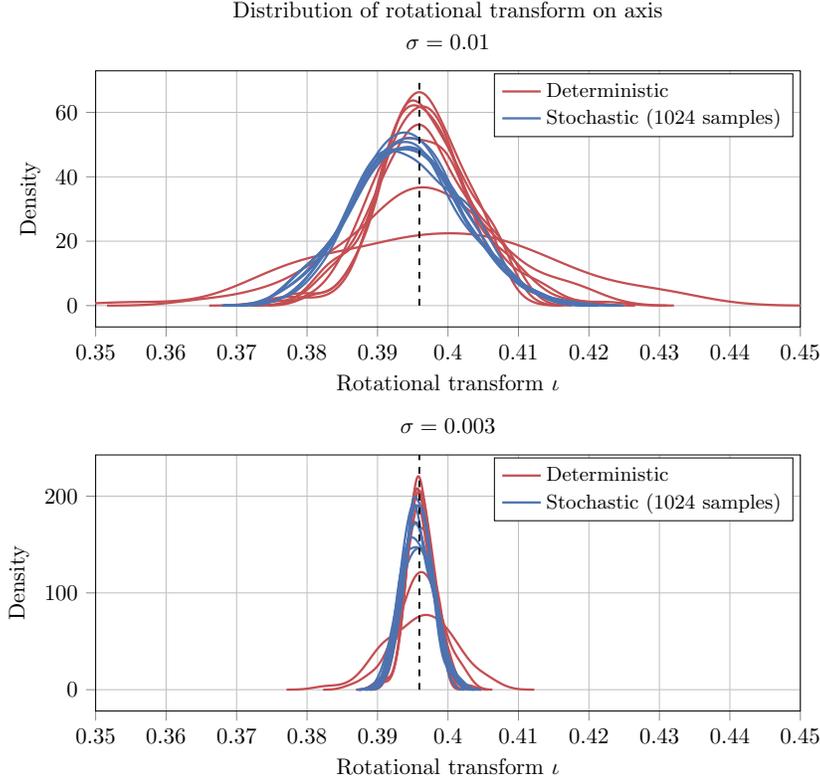

    \includestandalone[width=\textwidth]{./plots/iota_density}
    \caption{Distribution of the rotational transform $\iota$ on axis for each of the eight minimizers found from deterministic and stochastic optimization. The distribution is approximated using 128 independently drawn perturbed coil sets. The dashed line indicated the target rotational transform.}\label{fig:ncsx-iota-density}
\end{figure}

The objective includes a penalty that targets a certain rotational transform on the expansion axis.
To investigate the impact of coil perturbation errors on rotational transform, we draw 128 sets of perturbed coils, compute the magnetic axis for the resulting magnetic fields, and then compute the rotational transform $\iota$ on axis.
The resulting distribution of $\iota$ is shown in Figure~\ref{fig:ncsx-iota-density}.
As expected, for the perturbed coils the target rotational transform is not achieved exactly.
In agreement with our previous results, we also observe that the different minimizers obtained from deterministic optimization vary more than those obtained from stochastic optimization.

\subsection{Particle confinement}\label{subsec:conf}
\begin{figure}
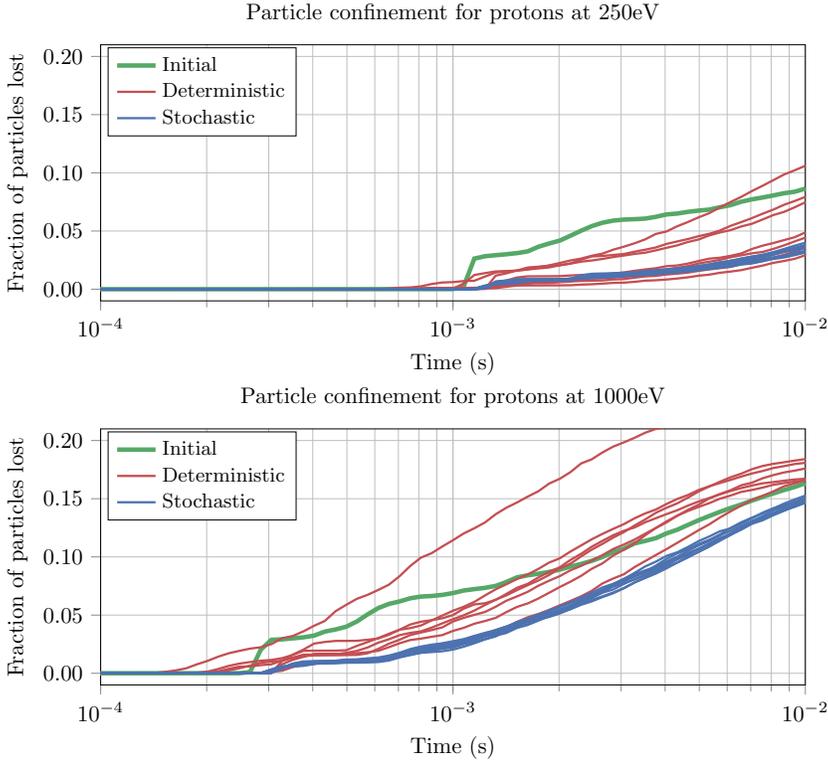

   \begin{center}
    \includestandalone[width=\textwidth]{./plots/confinement}
   \end{center} 
   \caption{Particle losses over time for 250eV and 1000eV protons spawned on axis, for perturbations of configurations obtained from stochastic and deterministic optimization.}\label{fig:confinement}
\end{figure}

As a final measure of performance, we compute particle trajectories for perturbations of the configurations obtained from deterministic and stochastic optimization. We consider both $250$eV protons and $1$keV protons. For reference, we note that a $3.4$keV proton in our designs has approximately
the same ratio of gyroradius to machine size as an energetic alpha particle in the ARIES-CS reactor.
We draw 10 perturbed coil configurations for each minimizer obtained from deterministic and stochastic optimization, as well as for the initial NCSX-like configuration.
We then spawn 1120 protons with random pitch angle on random points along the magnetic axis and follow them for 10ms, using the guiding center approximation~\cite{Boozer:1980,Littlejohn:1983,Cary:2009} without collisions.
Particles are considered lost if they move more than 30cm away from the axis.
Figure~\ref{fig:confinement} shows the average fraction of lost particles over time.

For both the lower energy protons and the higher energy protons considered here, we observe that the
configurations from stochastic optimization have better confinement than the initial configuration and the configurations obtained from deterministic optimization.  In addition, we see again that the different minimizers in the stochastic case all perform very
similarly, whereas those from deterministic optimization have highly
varying performance.  For the protons with lower energy, nearly all optimized configurations outperform the initial configuration. However, for the protons with
higher energy, this advantage is less clear. This suggests that quasi-symmetry only close to the axis is in general insufficient to
guarantee good particle confinement, and motivates ongoing work on direct coil optimization enforcing quasi-symmetry away from the axis.

\subsection{Risk-neutral vs.\ risk-averse optimization}\label{subsec:risk}

\begin{figure}[H]
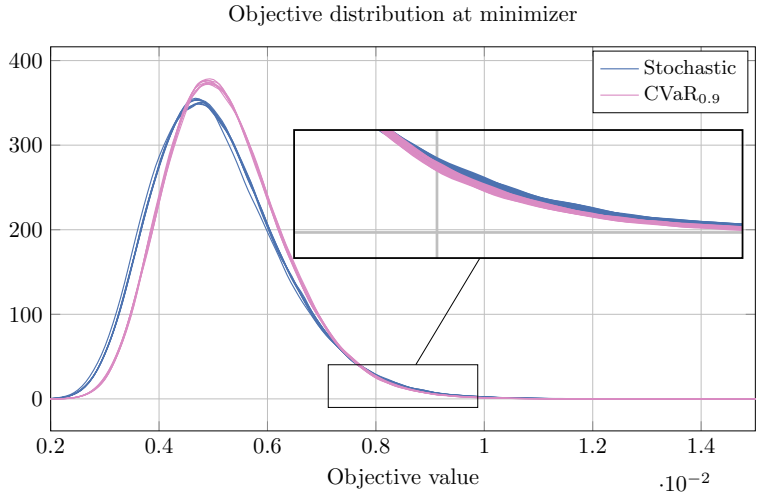

    \includestandalone[width=\textwidth]{./plots/riskaverse_density}
    \caption{Distribution of objective values for designs computed
      using risk-neutral stochastic formulation (blue) and CVaR
      risk-averse stochastic formulation (pink). $N_{\mathrm{MC}}=1024$ samples are
      used in the stochastic optimization, and the distributions are
      computed using 262\,144 independent samples. Slightly different
      designs are obtained using eight initializations for the
      optimization.}\label{fig:ncsx-riskaverse-dist}
\end{figure}

Finally, we compare the risk-neutral formulation (i.e.~minimization of
the expected value) with a risk-averse objective for an error
distribution with $\sigma=10^{-2}$.
We choose $\alpha=0.95$ and minimize $\mathrm{CVaR}_\alpha(f)$, meaning that we minimize the expected value of the tail containing the $5\%$ worst scenarios.
In Figure~\ref{fig:ncsx-riskaverse-dist} we show the out-of-sample distribution of the objective evaluated at the minimizers.
We can see that most of the distribution for the risk-neutral objective is closer to zero, while the tail is slightly thicker.
However, the difference is insignificant.
We attribute this to the quadratic penalty form of our objective: since all objective values are positive and the objective is squared, in order to control the mean large positive outliers have to be avoided.

\section{Conclusion and future work}\label{sec:conclusion}
We have extended the direct stellarator design approach of \cite{giuliani2020} to include random coil errors.
We emphasize that our formulation uses separate discretizations for the coils and their errors, which allows us to retain symmetries in the design space but considers perturbations that do not satisfy them.

We then studied and compared deterministic, risk-neutral, and risk-averse optimization for an NCSX-like example.
We found that the deterministic problem admits a large number of distinct minimizers which perform quite differently.
Including stochasticity reduces the number of different minimizers and results in minimizers that all perform nearly identical in terms of both objective value, as well as quasi-symmetry and rotational transform on and away from axis.
Moving from a risk-neutral to a risk-averse formulation does not result
in significantly different minimizing designs in our experiments: the tail of the distribution is reduced at the cost of somewhat worse average performance.
Finally, while we are able to achieve quasi-symmetry near the axis, this property is lost away from axis.
Current work is focused on including the non-quasi-symmetry measure presented in Figure~\ref{fig:ncsx-non-qs} on a range of surfaces to enforce quasi-symmetry away from axis.

\ack
This work was supported by a grant from the Simons Foundation (560651). AG is partially supported by a NSERC (Natural Sciences and Engineering Research Council of Canada) postdoctoral fellowship.
In addition, this work was supported in part through the NYU IT High Performance Computing resources, services, and staff expertise.

\section*{Code availability}
The code used to generate the numerical results is openly available under
\begin{center}
    \url{https://github.com/florianwechsung/PyPlasmaOpt/tree/fw/paper-stochastic/examples/stochastic} 
\end{center}

\printbibliography
\end{document}

%% file: latex-includes/makros.tex
\usepackage{bm}
\usepackage{xparse}
\usepackage{relsize}
\newcommand{\re}{\mathbb{R}}

\newcommand{\eps}{\varepsilon}
\newcommand{\epsb}{\boldsymbol{\varepsilon}}
\renewcommand{\l}{\left}
\renewcommand{\r}{\right}

\newcommand{\ab}{\mathbf{a}}

\newcommand{\cb}{\mathbf{c}}

\newcommand{\qb}{\mathbf{q}}

\newcommand{\xb}{\mathbf{x}}

\newcommand{\zb}{\mathbf{z}}

\newcommand{\Gammab}{\mathbf{\Gamma}}
\newcommand{\Xib}{\mathbf{\Xi}}
\newcommand{\zetab}{\boldsymbol{\zeta}}
\newcommand{\Bb}{\mathbf{B}}

\newcommand{\Sb}{\mathbf{S}}

\let\temp\phi
\let\phi\varphi
\let\varphi\temp

\newcommand{\vertiii}[1]{{\left\vert\kern-0.25ex\left\vert\kern-0.25ex\left\vert #1 
    \right\vert\kern-0.25ex\right\vert\kern-0.25ex\right\vert}}

\def\Proj_#1{\ensuremath{\operatorname{Proj}_{#1}\!}}

\newcommand{\dr}[1]{\ensuremath{\operatorname{d}\!{#1}}}

\def\drs_#1{\ensuremath{\operatorname{d}_{#1}\!}}
\def\Drs_#1{\ensuremath{\operatorname{D}_{#1}\!}}
\def\Drs^#1{\ensuremath{\operatorname{D}^{#1}\!}}

\def\Ers_#1{\mathrm{E}_{#1}}

\newcommand{\Cov}{\mathrm{Cov}}

\def\xunderbrace#1_#2{{\underbrace{#1}_{\mathclap{#2}}}}
\def\xoverbrace#1^#2{{\overbrace{#1}^{\mathclap{#2}}}}
\def\xunderarrow#1_#2{{\underset{\overset{\uparrow}{\mathclap{#2}}}{#1}}}
\def\xoverarrow#1^#2{{\overset{\underset{\downarrow}{\mathclap{#2}}}{#1}}}

\DeclareDocumentCommand{\boxedeq}{m o}{%
    \IfNoValueTF{#2}{%
        \rlap{\boxed{#1}}%
        \phantom{\hskip\fboxrule\hskip\fboxsep #1}%
        }{%
        \rlap{\boxed{#1#2}}%
        \phantom{\hskip\fboxrule\hskip\fboxsep #1}&\phantom{#2}%
    }%
}

\newcommand{\E}[1]{
    \mathbb{E}\l[#1\r]
}
\renewcommand{\P}[1]{
    \mathbb{P}\l(#1\r)
}
\newcommand{\cvar}[2]{
    \mathrm{CVaR}_{#1}\l[#2\r]
}

\newcommand{\pkg}[1]{\textsf{#1}}

%% file: plots/density.tex
\providecommand{\datadir}{.}%
\pgfplotstableread[col sep = semicolon, header=true]{\datadir/density.txt}\densitydata
\pgfplotsset{isstyle/.style={solid, line width=1.0pt}}
\pgfplotsset{oosstyle/.style={dashed, line width=1.0pt}}

\begin{tikzpicture}
    \begin{axis}[
      width=11cm,
      height=6cm,
      xlabel={Objective value},
      ylabel near ticks,
      ylabel={},
      ylabel style={overlay}, 
      yticklabel style={overlay},
      title={Objective distribution at minimizer},
      legend cell align=left,
      legend style={at={(0.99, 0.99)}, anchor=north east},
      legend columns=1,
      legend style={font=\small},
      scale only axis,
      grid=both,
      xmin=0.00,
      xmax=0.04,
      ]

      \addplot[isstyle, red] table[x=x_n_1_order_6_ig_0,y=y_n_1_order_6_ig_0] {\densitydata};
      \addlegendentry{Deterministic}
      \foreach \ig in {1,...,7} {
          \addplot[isstyle, red, forget plot] table[x=x_n_1_order_6_ig_\ig,y=y_n_1_order_6_ig_\ig] {\densitydata};
      }

      \addplot[isstyle, green] table[x=x_n_4_order_6_ig_0,y=y_n_4_order_6_ig_0] {\densitydata};
      \addlegendentry{Stochastic (4 samples)}
      \foreach \ig in {1,...,7} {
          \addplot[isstyle, green, forget plot] table[x=x_n_4_order_6_ig_\ig,y=y_n_4_order_6_ig_\ig] {\densitydata};
      }

      \addplot[isstyle, blue] table[x=x_n_1024_order_6_ig_0,y=y_n_1024_order_6_ig_0] {\densitydata};
      \addlegendentry{Stochastic (1024 samples)}
      \foreach \ig in {1,...,7} {
          \addplot[isstyle, blue, forget plot] table[x=x_n_1024_order_6_ig_\ig,y=y_n_1024_order_6_ig_\ig] {\densitydata};
      }

    \end{axis}
    \node at (13cm, 7cm)(a) {};
    \node at (-2cm, 0cm)(a) {};
\end{tikzpicture}

%% file: plots/qs_on_surface.tex
\providecommand{\datadir}{.}%

\pgfplotstableread[header=true]{\datadir/surface_data_0p01_0p01.txt}\qsdataone
\pgfplotstableread[header=true]{\datadir/surface_data_0p01_0p001.txt}\qsdataonesmaller
\pgfplotstableread[header=true]{\datadir/surface_data_0p003_0p003.txt}\qsdatathree
\pgfplotstableread[header=true]{\datadir/surface_data_0p003_0p0003.txt}\qsdatathreesmaller

\pgfplotsset{isstyle/.style={dashed, mark options={fill=none}, mark size=1pt, line width=1.0pt}}
\pgfplotsset{oosstyle/.style={solid, mark options={solid, fill=none}, mark size=1pt, line width=1.0pt}}

\begin{tikzpicture}
    \begin{groupplot}[
        group style={
            group name=my plots,
            group size=2 by 1,
            xlabels at=edge bottom,
            xticklabels at=edge bottom,
            horizontal sep=0pt
        },
        width=7.5cm,
        height=8cm,
        tickpos=left,
        xmin=-0.5, xmax=9.5,
        xlabel={Surface area},
        ylabel near ticks,
        legend cell align=left,
        legend style={at={(0.99, 0.01)}, anchor=south east},
        legend columns=1,
        legend style={font=\footnotesize},
        scale only axis,
        ymin=1e-4,
        ymax=3e-2,
        scale only axis,
        ylabel style={overlay}, yticklabel style={overlay}, xlabel style={overlay}, xticklabel style={overlay},
        trim axis left
        ]
    \nextgroupplot[ymode=log,grid=both,ytick align=outside, xtick align=outside, ylabel={$\|B_{\text{NotQS}}\|_{L^2}/\|B_{\text{QS}}\|_{L^2}$}]

      \addplot+[red, mark=*, oosstyle] table[x=target_areas,y=n_1_order_6_IG_0_oos_non_qs] {\qsdataone};
      \addlegendentry{Deterministic ($\sigma_\mathrm{opt}=0, \sigma_\mathrm{s} = 0.01$)}
      \foreach \ig in {1,...,7} {
          \addplot+[red, mark=*, oosstyle, forget plot] table[x=target_areas,y=n_1_order_6_IG_\ig_oos_non_qs] {\qsdataone};
      }
      \addplot+[red, mark=*, isstyle] table[x=target_areas,y=n_1_order_6_IG_0_oos_non_qs] {\qsdataonesmaller};
      \addlegendentry{Deterministic ($\sigma_\mathrm{opt}=0, \sigma_\mathrm{s} = 0.001$)}
      \foreach \ig in {1,...,7} {
          \addplot+[red, mark=*, isstyle, forget plot] table[x=target_areas,y=n_1_order_6_IG_\ig_oos_non_qs] {\qsdataonesmaller};
      }
      \addplot+[blue, mark=*, oosstyle] table[x=target_areas,y=n_1024_order_6_IG_0_oos_non_qs] {\qsdataone};
      \addlegendentry{Stochastic ($\sigma_\mathrm{opt}=0.01, \sigma_\mathrm{s} = 0.01$)}
      \foreach \ig in {1,...,7} {
          \addplot+[blue, mark=*, oosstyle, forget plot] table[x=target_areas,y=n_1024_order_6_IG_\ig_oos_non_qs] {\qsdataone};
      }
      \addplot+[blue, mark=*, isstyle] table[x=target_areas,y=n_1024_order_6_IG_0_oos_non_qs] {\qsdataonesmaller};
      \addlegendentry{Stochastic ($\sigma_\mathrm{opt}=0.01, \sigma_\mathrm{s} = 0.001$)}
      \foreach \ig in {1,...,7} {
          \addplot+[blue, mark=*, isstyle, forget plot] table[x=target_areas,y=n_1024_order_6_IG_\ig_oos_non_qs] {\qsdataonesmaller};
      }

      \nextgroupplot[ymode=log,grid=both,  ylabel near ticks, yticklabel pos=right]

      \addplot+[red, mark=*, oosstyle] table[x=target_areas,y=n_1_order_6_IG_0_oos_non_qs] {\qsdatathree};
      \addlegendentry{Deterministic ($\sigma_\mathrm{opt}=0, \sigma_\mathrm{s} = 0.003$)}
      \foreach \ig in {1,...,7} {
          \addplot+[red, mark=*, oosstyle, forget plot] table[x=target_areas,y=n_1_order_6_IG_\ig_oos_non_qs] {\qsdatathree};
      }
      \addplot+[isstyle, red, mark=*] table[x=target_areas,y=n_1_order_6_IG_0_oos_non_qs] {\qsdatathreesmaller};
      \addlegendentryexpanded{Deterministic ($\sigma_\mathrm{opt}=0, \sigma_\mathrm{s} = 0.0003$)}
      \foreach \ig in {1,...,7} {
          \addplot+[isstyle, red, mark=*, forget plot] table[x=target_areas,y=n_1_order_6_IG_\ig_oos_non_qs] {\qsdatathreesmaller};
      }
      \addplot+[blue, mark=*, oosstyle] table[x=target_areas,y=n_1024_order_6_IG_0_oos_non_qs] {\qsdatathree};
      \addlegendentry{Stochastic ($\sigma_\mathrm{opt}=0.003, \sigma_\mathrm{s} = 0.003$)}
      \foreach \ig in {1,...,7} {
          \addplot+[blue, mark=*, oosstyle, forget plot] table[x=target_areas,y=n_1024_order_6_IG_\ig_oos_non_qs] {\qsdatathree};
      }
      \addplot+[isstyle, blue, mark=*] table[x=target_areas,y=n_1024_order_6_IG_0_oos_non_qs] {\qsdatathreesmaller};
      \addlegendentryexpanded{Stochastic ($\sigma_\mathrm{opt}=0.003, \sigma_\mathrm{s} = 0.0003$)}
      \foreach \ig in {1,...,7} {
          \addplot+[isstyle, blue, mark=*, forget plot] table[x=target_areas,y=n_1024_order_6_IG_\ig_oos_non_qs] {\qsdatathreesmaller};
      }
    \end{groupplot}
    \node (title) at (7.5cm, 8.0cm) [above, yshift=\pgfkeysvalueof{/pgfplots/every axis title shift}] {Non quasi-symmetric part of $B$ away from axis};
    \node at (16.5cm, 7cm)(a) {};
    \node at (-1.5cm, -1cm)(a) {};
\end{tikzpicture}

%% file: plots/iota_density.tex
\providecommand{\datadir}{.}%
\pgfplotstableread[col sep = semicolon, header=true]{\datadir/iota_density_0p01.txt}\densitydata
\pgfplotstableread[col sep = semicolon, header=true]{\datadir/iota_density_0p003.txt}\densitydatasmall
\pgfplotsset{isstyle/.style={solid, line width=1.0pt}}
\pgfplotsset{oosstyle/.style={dashed, line width=1.0pt}}

\begin{tikzpicture}
    \begin{groupplot}[
        group style={
            group name=my plots,
            group size=1 by 2,
            vertical sep=2.0cm
        },
        xmin=0.35,
        xmax=0.45,
        width=11cm,
        height=4cm,
        tickpos=left,
        ylabel near ticks,
        legend cell align=left,
        legend style={at={(0.99, 0.99)}, anchor=north east},
        legend columns=1,
        legend style={font=\small},
        scale only axis,
        scale only axis,
        ylabel style={overlay}, yticklabel style={overlay}, 
        trim axis left,
        xlabel={Rotational transform $\iota$},
        ]
      \nextgroupplot[grid=both,ytick align=outside, xtick align=outside, ylabel={Density}, title={$\sigma=0.01$}]
      \addplot[isstyle, red] table[x=x_n_1_order_6_ig_0,y=y_n_1_order_6_ig_0] {\densitydata};
      \addlegendentry{Deterministic}
      \foreach \ig in {1,...,7} {
          \addplot[isstyle, red, forget plot] table[x=x_n_1_order_6_ig_\ig,y=y_n_1_order_6_ig_\ig] {\densitydata};
      }

      \addplot[isstyle, blue] table[x=x_n_1024_order_6_ig_0,y=y_n_1024_order_6_ig_0] {\densitydata};
      \addlegendentry{Stochastic (1024 samples)}
      \foreach \ig in {1,...,7} {
          \addplot[isstyle, blue, forget plot] table[x=x_n_1024_order_6_ig_\ig,y=y_n_1024_order_6_ig_\ig] {\densitydata};
      }
      \draw [thick, dashed] (axis cs:0.395938929522566,0) -- (axis cs:0.395938929522566,70);

      \nextgroupplot[grid=both,ytick align=outside, xtick align=outside, ylabel={Density}, title={$\sigma=0.003$}]
      \addplot[isstyle, red] table[x=x_n_1_order_6_ig_0,y=y_n_1_order_6_ig_0] {\densitydatasmall};
      \addlegendentry{Deterministic}
      \foreach \ig in {1,...,7} {
          \addplot[isstyle, red, forget plot] table[x=x_n_1_order_6_ig_\ig,y=y_n_1_order_6_ig_\ig] {\densitydatasmall};
      }

      \addplot[isstyle, blue] table[x=x_n_1024_order_6_ig_0,y=y_n_1024_order_6_ig_0] {\densitydatasmall};
      \addlegendentry{Stochastic (1024 samples)}
      \foreach \ig in {1,...,7} {
          \addplot[isstyle, blue, forget plot] table[x=x_n_1024_order_6_ig_\ig,y=y_n_1024_order_6_ig_\ig] {\densitydatasmall};
      }
        \draw [thick, dashed] (axis cs:0.395938929522566,0) -- (axis cs:0.395938929522566,250);
\end{groupplot}
    \node (title) at (5.5cm, 4.5cm) [above, yshift=\pgfkeysvalueof{/pgfplots/every axis title shift}] {Distribution of rotational transform on axis};
      



    \node at (13cm, 5cm)(a) {};
    \node at (-2cm, 0cm)(a) {};
\end{tikzpicture}

%% file: plots/confinement.tex
\providecommand{\datadir}{.}%
\pgfplotstableread[col sep = semicolon, header=true]{\datadir/confinement_250.txt}\confinementtwofifty
\pgfplotstableread[col sep = semicolon, header=true]{\datadir/confinement_1000.txt}\confinementthousand
\pgfplotsset{isstyle/.style={solid, line width=1.0pt}}
\pgfplotsset{oosstyle/.style={dashed, line width=1.0pt}}

\begin{tikzpicture}
    \begin{groupplot}[
        group style={
            group name=my plots,
            group size=1 by 2,
            vertical sep=2.0cm
        },
        xmin=1e-4,
        xmax=1e-2,
        ymin=-0.01,
        ymax=0.21,
        width=11cm,
        height=4cm,
        tickpos=left,
        ylabel near ticks,
        legend cell align=left,
        legend style={at={(0.01, 0.99)}, anchor=north west},
        legend columns=1,
        legend style={font=\small},
        scale only axis,
        scale only axis,
        ylabel style={overlay},
        yticklabel style={overlay,
            /pgf/number format/fixed,
            /pgf/number format/precision=2,
            /pgf/number format/fixed zerofill
        },
        trim axis left,
        xlabel={Time (s)},
        ]
      \nextgroupplot[xmode=log, grid=both,ytick align=outside, xtick align=outside, ylabel={Fraction of particles lost}, title={Particle confinement for protons at 250eV}]
      \addplot[isstyle, green, line width=2.0] table[x=t,y=NCSX_Initial_(Perturbed)] {\confinementtwofifty};
      \addlegendentry{Initial}
      \addplot[isstyle, red] table[x=t,y=Deterministic_(ig_0)_Perturbed] {\confinementtwofifty};
      \addlegendentry{Deterministic}
      \foreach \ig in {1,...,7} {
          \addplot[isstyle, red, forget plot] table[x=t,y=Deterministic_(ig_\ig)_Perturbed] {\confinementtwofifty};
      }
      \addplot[isstyle, blue] table[x=t,y=Stochastic_(ig_0)_Perturbed] {\confinementtwofifty};
      \addlegendentry{Stochastic}
      \foreach \ig in {1,...,7} {
          \addplot[isstyle, blue, forget plot] table[x=t,y=Stochastic_(ig_\ig)_Perturbed] {\confinementtwofifty};
      }
      \nextgroupplot[xmode=log, grid=both,ytick align=outside, xtick align=outside, ylabel={Fraction of particles lost}, title={Particle confinement for protons at 1000eV}]
      \addplot[isstyle, green, line width=2.0] table[x=t,y=NCSX_Initial_(Perturbed)] {\confinementthousand};
      \addlegendentry{Initial}
      \addplot[isstyle, red] table[x=t,y=Deterministic_(ig_0)_Perturbed] {\confinementthousand};
      \addlegendentry{Deterministic}
      \foreach \ig in {1,...,6} {
          \addplot[isstyle, red, forget plot] table[x=t,y=Deterministic_(ig_\ig)_Perturbed] {\confinementthousand};
      }
      \addplot[isstyle, blue] table[x=t,y=Stochastic_(ig_0)_Perturbed] {\confinementthousand};
      \addlegendentry{Stochastic}
      \foreach \ig in {1,...,7} {
          \addplot[isstyle, blue, forget plot] table[x=t,y=Stochastic_(ig_\ig)_Perturbed] {\confinementthousand};
      }




\end{groupplot}
    \node at (13cm, 5cm)(a) {};
    \node at (-2cm, 0cm)(a) {};
\end{tikzpicture}

%% file: plots/riskaverse_density.tex
\providecommand{\datadir}{.}%
\pgfplotstableread[col sep = semicolon, header=true]{\datadir/riskaverse_density.txt}\densitydata
\pgfplotsset{isstyle/.style={solid, line width=0.5pt}}
\pgfplotsset{oosstyle/.style={dashed, line width=0.5pt}}

\begin{tikzpicture}[spy using outlines={rectangle, magnification=3, width=7cm, height=2cm, connect spies}]
    \begin{axis}[
      width=11cm,
      height=6cm,
      xlabel={Objective value},
      ylabel near ticks,
      ylabel={},
      ylabel style={overlay}, 
      yticklabel style={overlay},
      title={Objective distribution at minimizer},
      legend cell align=left,
      legend style={at={(0.99, 0.99)}, anchor=north east},
      legend columns=1,
      legend style={font=\small},
      scale only axis,
      grid=both,
      xmin=0.002,
      xmax=0.015,
      ]

      \addplot[isstyle, blue] table[x=x_mode_stochastic_ig_0,y=y_mode_stochastic_ig_0] {\densitydata};
      \addlegendentry{Stochastic}
      \foreach \ig in {1,...,7} {
          \addplot[isstyle, blue, forget plot] table[x=x_mode_stochastic_ig_\ig,y=y_mode_stochastic_ig_\ig] {\densitydata};
      }

      \addplot[isstyle, pink] table[x=x_mode_cvar0p9_ig_0,y=y_mode_cvar0p9_ig_0] {\densitydata};
      \addlegendentry{$\mathrm{CVaR}_{0.9}$}
      \foreach \ig in {1,...,7} {
          \addplot[isstyle, pink, forget plot] table[x=x_mode_cvar0p9_ig_\ig,y=y_mode_cvar0p9_ig_\ig] {\densitydata};
      }

    \end{axis}
    \spy on (5.5, 0.7) in node [] at (7.3, 3.7);
    \node at (13cm, 7cm)(a) {};
    \node at (-2cm, 0cm)(a) {};
\end{tikzpicture}

%% file: refs.bib
@article{Pedersen_2006,
	title = {Experimental demonstration of a compact stellarator magnetic trap using four circular coils},
	volume = {13},
	issn = {1070-664X, 1089-7674},
	doi = {10.1063/1.2149313},
	number = {1},
	journal = {Physics of Plasmas},
	author = {Pedersen, T. Sunn and Kremer, J. P. and Lefrancois, R. G. and Marksteiner, Q. and Sarasola, X. and Ahmad, N.},
	month = jan,
	year = {2006},
}

@article{Rummel_2004,
	title = {Accuracy of the {Construction} of the {Superconducting} {Coils} for {WENDELSTEIN} 7-{X}},
	volume = {14},
	issn = {1051-8223},
	doi = {10.1109/TASC.2004.830584},
	language = {en},
	number = {2},
	journal = {IEEE Transactions on Appiled Superconductivity},
	author = {Rummel, T. and Risse, K. and Viebke, H. and Braeuer, T. and Kisslinger, J.},
	month = jun,
	year = {2004},
	pages = {1394--1398},
}

@inproceedings{Kremer_2003,
	title = {The {Status} of the {Design} and {Construction} of the {Columbia} {Non}-neutral {Torus}},
	volume = {692},
	booktitle = {{AIP} {Conference} {Proceedings}},
	publisher = {American Institute of Physics},
	author = {Kremer, JP and Pedersen, Thomas Sunn and Pomphrey, N and Reiersen, W and Dahlgren, F},
	year = {2003},
	note = {Issue: 1},
	pages = {320--325},
}

@article{virtanen2020scipy,
  title={SciPy 1.0: fundamental algorithms for scientific computing in Python},
  author={Virtanen, Pauli and Gommers, Ralf and Oliphant, Travis E and Haberland, Matt and Reddy, Tyler and Cournapeau, David and Burovski, Evgeni and Peterson, Pearu and Weckesser, Warren and Bright, Jonathan and others},
  journal={Nature Methods},
  volume={17},
  number={3},
  pages={261--272},
  year={2020},
  publisher={Nature Publishing Group}
}

@article{strickler_designing_2002,
	title = {Designing {Coils} for {Compact} {Stellarators}},
	volume = {41},
	doi = {10.13182/FST02-A206},
	number = {2},
	journal = {Fusion Science and Technology},
	author = {Strickler, Dennis J. and Berry, Lee A. and Hirshman, Steven P.},
	month = mar,
	year = {2002},
	pages = {107--115},
}

@inproceedings{drevlak_optimization_1999,
	title = {Optimization of heterogenous magnet systems},
	booktitle = {Proceedings of the 12th {International} {Stellarator} {Workshop}},
	author = {Drevlak, Michael},
	year = {1999},
}

@article{zhu2017,
	title = {New method to design stellarator coils without the winding surface},
	volume = {58},
	issn = {0029-5515, 1741-4326},
	doi = {10.1088/1741-4326/aa8e0a},
	number = {1},
	journal = {Nuclear Fusion},
	author = {Zhu, Caoxiang and Hudson, Stuart R. and Song, Yuntao and Wan, Yuanxi},
	month = jan,
	year = {2018},
	pages = {016008},
}

@article{Zhu_2018,
	doi = {10.1088/1361-6587/aab8c2},
	year = 2018,
	month = {apr},
	publisher = {{IOP} Publishing},
	volume = {60},
	number = {6},
	pages = {065008},
	author = {Caoxiang Zhu and Stuart R Hudson and Yuntao Song and Yuanxi Wan},
	title = {Designing stellarator coils by a modified Newton method using {FOCUS}},
	journal = {Plasma Physics and Controlled Fusion},
}

@article{bader_2019, title={Stellarator equilibria with reactor relevant energetic particle losses}, volume={85}, DOI={10.1017/S0022377819000680}, number={5}, journal={Journal of Plasma Physics}, publisher={Cambridge University Press}, author={Bader, Aaron and Drevlak, M. and Anderson, D. T. and Faber, B. J. and Hegna, C. C. and Likin, K. M. and Schmitt, J. C. and Talmadge, J. N.}, year={2019}, pages={905850508}}

@article{gates_recent_2017,
	title = {Recent advances in stellarator optimization},
	volume = {57},
	doi = {10.1088/1741-4326/aa8ba0},
	number = {12},
	journal = {Nuclear Fusion},
	author = {Gates, D.A. and Boozer, A.H. and Brown, T. and Breslau, J. and Curreli, D. and Landreman, M. and Lazerson, S.A. and Lore, J. and Mynick, H. and Neilson, G.H. and Pomphrey, N. and Xanthopoulos, P. and Zolfaghari, A.},
	month = dec,
	year = {2017},
	pages = {126064},
}

@article{gates2018stellarator,
  title={Stellarator research opportunities: a report of the National Stellarator Coordinating Committee},
  author={Gates, David A and Anderson, David and Anderson, Seth and Zarnstorff, M and Spong, Donald A and Weitzner, Harold and Neilson, GH and Ruzic, D and Andruczyk, D and Harris, Jeffrey H and others},
  journal={Journal of Fusion Energy},
  volume={37},
  number={1},
  pages={51--94},
  year={2018},
  publisher={Springer}
}

@article{chen_new_2011,
	title = {A new level-set based approach to shape and topology optimization under geometric uncertainty},
	volume = {44},
	issn = {1615-147X, 1615-1488},
	doi = {10.1007/s00158-011-0660-9},
	number = {1},
	journal = {Structural and Multidisciplinary Optimization},
	author = {Chen, Shikui and Chen, Wei},
	month = jul,
	year = {2011},
	pages = {1--18},
}

@inproceedings{wang_conditional_2011,
	address = {Orlando, Florida},
	title = {Conditional sampling and experiment design for quantifying manufacturing error of transonic airfoil},
	isbn = {978-1-60086-950-1},
	doi = {10.2514/6.2011-658},
	booktitle = {49th {AIAA} {Aerospace} {Sciences} {Meeting} including the {New} {Horizons} {Forum} and {Aerospace} {Exposition}},
	publisher = {American Institute of Aeronautics and Astronautics},
	author = {Wang, Qiqi and Chen, Han and Hu, Rui and Constantine, Paul},
	month = jan,
	year = {2011},
}

@article{liu_quantification_2017,
	title = {Quantification of {Airfoil} {Geometry}-{Induced} {Aerodynamic} {Uncertainties}---{Comparison} of {Approaches}},
	volume = {5},
	issn = {2166-2525},
	doi = {10.1137/15M1050239},
	number = {1},
	journal = {SIAM/ASA Journal on Uncertainty Quantification},
	author = {Liu, Dishi and Litvinenko, Alexander and Schillings, Claudia and Schulz, Volker},
	month = jan,
	year = {2017},
	pages = {334--352},
}

@article{lobsien_physics_2020,
	title = {Physics analysis of results of stochastic and classic stellarator coil optimization},
	volume = {60},
	doi = {10.1088/1741-4326/ab7211},
	number = {4},
	journal = {Nuclear Fusion},
	author = {Lobsien, Jim-Felix and Drevlak, Michael and Jenko, Frank and Maurer, Maurice and Navarro, Alejandro Bañon and Nührenberg, Carolin and Pedersen, Thomas Sunn and Smith, Håkan M. and Turkin, Yuriy and {W7-X Team}},
	month = apr,
	year = {2020},
}

@article{lobsien_stellarator_2018,
	title = {Stellarator coil optimization towards higher engineering tolerances},
	volume = {58},
	doi = {10.1088/1741-4326/aad431},
	number = {10},
	journal = {Nuclear Fusion},
	author = {Lobsien, Jim-Felix and Drevlak, Michael and Sunn Pedersen, Thomas and {W7-X Team}},
	month = oct,
	year = {2018},
}

@article{lobsien_improved_2020,
	title = {Improved performance of stellarator coil design optimization},
	volume = {86},
	doi = {10.1017/S0022377820000227},
	number = {2},
	journal = {Journal of Plasma Physics},
	author = {Lobsien, Jim-Felix and Drevlak, Michael and Kruger, Thomas and Lazerson, Samuel and Zhu, Caoxiang and Pedersen, Thomas Sunn},
	month = apr,
	year = {2020},
}

@book{adler_geometry_2010,
	title = {The {Geometry} of {Random} {Fields}},
	language = {en},
	publisher = {Society for Industrial and Applied Mathematics},
	author = {Adler, Robert J},
	month = jan,
	year = {2010},
	doi = {10.1137/1.9780898718980},
}

@book{rasmussen_gaussian_2006,
	address = {Cambridge, Mass},
	series = {Adaptive computation and machine learning},
	title = {Gaussian processes for machine learning},
	publisher = {MIT Press},
	author = {Rasmussen, Carl Edward and Williams, Christopher K. I.},
	year = {2006},
}

@book{scholkopf_learning_2002,
	address = {Cambridge, Mass},
	series = {Adaptive computation and machine learning},
	title = {Learning with kernels: support vector machines, regularization, optimization, and beyond},
	shorttitle = {Learning with kernels},
	publisher = {MIT Press},
	author = {Schölkopf, Bernhard and Smola, Alexander J.},
	year = {2002},
}

@article{rockafellar_optimization_2000,
	title = {Optimization of conditional value-at-risk},
	volume = {2},
	issn = {14651211},
	doi = {10.21314/JOR.2000.038},
	language = {en},
	number = {3},
	journal = {The Journal of Risk},
	author = {Rockafellar, R. Tyrrell and Uryasev, Stanislav},
	year = {2000},
	pages = {21--41},
}

@article{giuliani2020,
  title={Single-stage gradient-based stellarator coil design: Optimization for near-axis quasi-symmetry},
  author={Giuliani, Andrew and Wechsung, Florian and Cerfon, Antoine and Stadler, Georg and Landreman, Matt},
  journal={arXiv preprint arXiv:2010.02033},
  year={2020}
}

@article{Boozer:1980,
author = {Boozer,Allen H. },
title = {Guiding center drift equations},
journal = {The Physics of Fluids},
volume = {23},
number = {5},
pages = {904-908},
year = {1980},
doi = {10.1063/1.863080},

URL = { 
        https://aip.scitation.org/doi/abs/10.1063/1.863080
    
},
eprint = { 
        https://aip.scitation.org/doi/pdf/10.1063/1.863080
    
}

}

@article{Landreman_2018_sensitivity,
	doi = {10.1088/1741-4326/aac197},
	url = {https://doi.org/10.1088/1741-4326/aac197},
	year = 2018,
	month = {jun},
	publisher = {{IOP} Publishing},
	volume = {58},
	number = {7},
	pages = {076023},
	author = {Matt Landreman and Elizabeth Paul},
	title = {Computing local sensitivity and tolerances for stellarator physics properties using shape gradients},
	journal = {Nuclear Fusion}
}

@article{Littlejohn:1983, title={Variational principles of guiding centre motion}, volume={29}, DOI={10.1017/S002237780000060X}, number={1}, journal={Journal of Plasma Physics}, publisher={Cambridge University Press}, author={Littlejohn, R. G.}, year={1983}, pages={111–125}}

@article{Cary:2009,
  title = {Hamiltonian theory of guiding-center motion},
  author = {Cary, John R. and Brizard, Alain J.},
  journal = {Rev. Mod. Phys.},
  volume = {81},
  issue = {2},
  pages = {693--738},
  numpages = {0},
  year = {2009},
  month = {May},
  publisher = {American Physical Society},
  doi = {10.1103/RevModPhys.81.693},
  url = {https://link.aps.org/doi/10.1103/RevModPhys.81.693}
}

@article{Carlton2021,
  title={Computing the shape gradient of stellarator coil complexity with respect to the plasma boundary},
  author={Carlton-Jones, Arthur and Paul, Elizabeth J and Dorland, William},
  journal={Journal of Plasma Physics},
  volume={87},
  number={2},
  year={2021},
  publisher={Cambridge University Press}
}

@article{Drevlak_2013,
author = {Drevlak, M. and Brochard, F. and Helander, P. and Kisslinger, J. and Mikhailov, M. and N\"uhrenberg, C. and Nührenberg, J. and Turkin, Y.},
title = {ESTELL: A Quasi-Toroidally Symmetric Stellarator},
journal = {Contributions to Plasma Physics},
volume = {53},
number = {6},
pages = {459-468},
doi = {10.1002/ctpp.201200055},
year = {2013}
}

@article{geraldini_2021, title={An adjoint method for determining the sensitivity of island size to magnetic field variations}, volume={87}, DOI={10.1017/S0022377821000428}, number={3}, journal={Journal of Plasma Physics}, publisher={Cambridge University Press}, author={Geraldini, Alessandro and Landreman, M. and Paul, E.}, year={2021}, pages={905870302}}

@article{Henneberg_2019,
	doi = {10.1088/1741-4326/aaf604},
	year = 2019,
	publisher = {{IOP} Publishing},
	volume = {59},
	number = {2},
	pages = {026014},
	author = {S.A. Henneberg and M. Drevlak and C. N{\"u}hrenberg and C.D. Beidler and Y. Turkin and J. Loizu and P. Helander},
	title = {Properties of a new quasi-axisymmetric configuration},
	journal = {Nuclear Fusion}
}

@article{Hudson2018,
title = {Differentiating the shape of stellarator coils with respect to the plasma boundary},
journal = {Physics Letters A},
volume = {382},
number = {38},
pages = {2732-2737},
year = {2018},
issn = {0375-9601},
doi = {https://doi.org/10.1016/j.physleta.2018.07.016},
author = {S.R. Hudson and C. Zhu and D. Pfefferlé and L. Gunderson},
}

@article{Landreman2016,
author = {Landreman,Matt  and Boozer,Allen H. },
title = {Efficient magnetic fields for supporting toroidal plasmas},
journal = {Physics of Plasmas},
volume = {23},
number = {3},
pages = {032506},
year = {2016},
doi = {10.1063/1.4943201},
}

@article{Landreman_2017,
	doi = {10.1088/1741-4326/aa57d4},
	year = 2017,
	month = {feb},
	publisher = {{IOP} Publishing},
	volume = {57},
	number = {4},
	pages = {046003},
	author = {Matt Landreman},
	title = {An improved current potential method for fast computation of stellarator coil shapes},
	journal = {Nuclear Fusion}
}

@article{landreman_2018, title={Direct construction of optimized stellarator shapes. {P}art 1. {T}heory in cylindrical coordinates}, volume={84}, DOI={10.1017/S0022377818001289}, number={6}, journal={Journal of Plasma Physics}, publisher={Cambridge University Press}, author={Landreman, Matt and Sengupta, Wrick}, year={2018}, pages={905840616}}

@article{landreman_2019, title={Direct construction of optimized stellarator shapes. {P}art 2. {N}umerical quasisymmetric solutions}, volume={85}, DOI={10.1017/S0022377818001344}, number={1}, journal={Journal of Plasma Physics}, publisher={Cambridge University Press}, author={Landreman, Matt and Sengupta, Wrick and Plunk, Gabriel G.}, year={2019}, pages={905850103}}

@article{Paul_2018,
	doi = {10.1088/1741-4326/aac1c7},
	year = 2018,
	month = {may},
	publisher = {{IOP} Publishing},
	volume = {58},
	number = {7},
	pages = {076015},
	author = {E.J. Paul and M. Landreman and A. Bader and W. Dorland},
	title = {An adjoint method for gradient-based optimization of stellarator coil shapes},
	journal = {Nuclear Fusion}
}

@article{Zhu_2018_error,
	doi = {10.1088/1361-6587/aab6cb},
	year = 2018,
	month = {apr},
	publisher = {{IOP} Publishing},
	volume = {60},
	number = {5},
	pages = {054016},
	author = {Caoxiang Zhu and Stuart R Hudson and Samuel A Lazerson and Yuntao Song and Yuanxi Wan},
	title = {Hessian matrix approach for determining error field sensitivity to coil deviations},
	journal = {Plasma Physics and Controlled Fusion}
}

@article{Zhu_2019,
	doi = {10.1088/1741-4326/ab3a7c},
	year = 2019,
	month = {sep},
	publisher = {{IOP} Publishing},
	volume = {59},
	number = {12},
	pages = {126007},
	author = {Caoxiang Zhu and David A. Gates and Stuart R. Hudson and Haifeng Liu and Yuhong Xu and Akihiro Shimizu and Shoichi Okamura},
	title = {Identification of important error fields in stellarators using the Hessian matrix method},
	journal = {Nuclear Fusion}
}

@article{Klinger2013,
title = "Towards assembly completion and preparation of experimental campaigns of {W}endelstein {7-X} in the perspective of a path to a stellarator fusion power plant",
journal = "Fusion Engineering and Design",
volume = "88",
number = "6",
pages = "461 - 465",
year = "2013",
note = "Proceedings of the 27th Symposium On Fusion Technology (SOFT-27); Liège, Belgium, September 24-28, 2012",
issn = "0920-3796",
doi = "https://doi.org/10.1016/j.fusengdes.2013.02.153",
author = "T. Klinger and C. Baylard and C.D. Beidler and J. Boscary and H.S. Bosch and A. Dinklage and D. Hartmann and P. Helander and H. Ma{\ss}berg and A. Peacock and T.S. Pedersen and T. Rummel and F. Schauer and L. Wegener and R. Wolf"
}

@INPROCEEDINGS{Strykowsky09,

  author={R. L. {Strykowsky} and T. {Brown} and J. {Chrzanowski} and M. {Cole} and P. {Heitzenroeder} and G. H. {Neilson} and D. {Rej} and M. {Viol}},

  booktitle={2009 23rd IEEE/NPSS Symposium on Fusion Engineering}, 

  title={Engineering cost   schedule lessons learned on {NCSX}}, 

  year={2009},

  volume={},

  number={},

  pages={1-4},}

@ARTICLE{Neilson10,

  author={G. H. {Neilson} and C. O. {Gruber} and J. H. {Harris} and D. J. {Rej} and R. T. {Simmons} and R. L. {Strykowsky}},

  journal={IEEE Transactions on Plasma Science}, 

  title={Lessons Learned in Risk Management on {NCSX}}, 

  year={2010},

  volume={38},

  number={3},

  pages={320-327},}

@article{Dewar1998,
title = "Stellarator symmetry",
journal = "Physica D: Nonlinear Phenomena",
volume = "112",
number = "1",
pages = "275--280",
year = "1998",
note = "Proceedings of the Workshop on Time-Reversal Symmetry in Dynamical Systems",
issn = "0167-2789",
doi = "https://doi.org/10.1016/S0167-2789(97)00216-9",
author = "R.L. Dewar and S.R. Hudson",
}

@book{RasmussenWilliams06,
  title = {Gaussian Processes for Machine Learning},
  author = {Rasmussen, CE. and Williams, CKI.},
  pages = {248},
  series = {Adaptive Computation and Machine Learning},
  publisher = {MIT Press},
  organization = {Max-Planck-Gesellschaft},
  school = {Biologische Kybernetik},
  address = {Cambridge, MA, USA},
  month = jan,
  year = {2006},
  month_numeric = {1}
}

@Book{ShapiroDentchevaRuszczynski09,
  Title                    = {Lectures on Stochastic Programming: Modeling and Theory},
  Author                   = {Alexander Shapiro and Darinka Dentcheva and Andrezj Ruszczynski},
  Publisher                = {Society for Industrial and Applied Mathematics},
  Year                     = {2009}
}

@article{KrokhmalPalmquistUryasev02,
  title={Portfolio optimization with conditional value-at-risk objective and constraints},
  author={Krokhmal, Pavlo and Palmquist, Jonas and Uryasev, Stanislav},
  journal={Journal of Risk},
  volume={4},
  pages={43--68},
  year={2002},
  publisher={Citeseer}
}

@Article{KouriSurowiec16,
  Title                    = {Risk-Averse {PDE}-Constrained Optimization Using the Conditional Value-At-Risk},
  Author                   = {D. P. Kouri and T. M. Surowiec},
  Journal                  = {SIAM Journal on Optimization},
  Year                     = {2016},
  Number                   = {1},
  Pages                    = {365-396},
  Volume                   = {26},

  Doi                      = {10.1137/140954556}
}
